\makeatletter \let\cl@chapter\relax \makeatother

\documentclass[numbook]{svjour3}

\journalname{Numerical Algorithms}

\usepackage[numbers,sort&compress]{natbib}
\usepackage{mathptmx}

\usepackage{multicol}

\usepackage{microtype}

\usepackage{appendix}

\usepackage{amsmath}
\usepackage{amsfonts}
\usepackage{amssymb}

\allowdisplaybreaks[1]

\usepackage{hyperref}

\numberwithin{equation}{section}

\usepackage[capitalize]{cleveref}

\providecommand{\Ent}[1]{\lfloor #1 \rfloor}

\begin{document}

\thispagestyle{empty}

\title{Construction of New Generalizations of Wynn's Epsilon and Rho
  Algorithm by Solving Finite Difference Equations in the Transformation
  Order}

\titlerunning{New Generalizations of Wynn's Epsilon and Rho Algorithm}

\author{Xiang-Ke Chang \and
        Yi He \and
        Xing-Biao Hu \and
        Jian-Qing Sun \and
        Ernst Joachim Weniger}

\institute{Xiang-Ke Chang \at
  LSEC, ICMSEC, Academy of Mathematics and Systems Science, Chinese
  Academy of Sciences, P.O.Box 2719, Beijing 100190, PR China;
  and School of Mathematical Sciences, University of Chinese Academy of
  Sciences, Beijing 100049, PR China \\
  \email{changxk@lsec.cc.ac.cn}
  \and
  Yi He \at
  Wuhan Institute of Physics and Mathematics, Chinese Academy of
  Sciences, Wuhan 430071, PR China \\
  \email{heyi@wipm.ac.cn}
  \and
   Xing-Biao Hu \at
  LSEC, ICMSEC, Academy of Mathematics and Systems Science, Chinese
  Academy of Sciences, P.O.Box 2719, Beijing 100190, PR China; and School
  of Mathematical Sciences, University of Chinese Academy of Sciences,
  Beijing 100049, PR China \\
  \email{hxb@lsec.cc.ac.cn}
  \and
  Jian-Qing Sun \at
  School of Mathematical Sciences, Ocean University of China, Qingdao,
  266100, PR China \\
  \email{sunjq@lsec.cc.ac.cn}
  \and
  Ernst Joachim Weniger \at
  Instit{u}t f\"{u}r Physikalische und Theoretische Chemie,
  Universit\"{a}t Regensburg, D-93040 Regensburg, Germany \\
  \email{joachim.weniger@chemie.uni-regensburg.de}}

\bigskip \bigskip

\date{To Appear in Numerical Algorithms: \today}

\maketitle

\begin{abstract}
  We construct new sequence transformations based on Wynn's epsilon and
  rho algorithms. The recursions of the new algorithms include the
  recursions of Wynn's epsilon and rho algorithm and of Osada's
  generalized rho algorithm as special cases. We demonstrate the
  performance of our algorithms numerically by applying them to some
  linearly and logarithmically convergent sequences as well as some
  divergent series.
\end{abstract}



\vspace*{0.5cm}

\begin{center}
  \textbf{We dedicate this article to the memory of Peter Wynn (1931 -
    2017)}
\end{center}

 \newpage

\tableofcontents

\newpage


%
\typeout{==> Section: Introduction}
\section{Introduction}
\label{Sec:Intro}
%

Sequences and series are extremely important mathematical tools. They
appear naturally in many methods for solving differential and integral
equations, in discretization methods, in quadrature schemes, in
perturbation techniques, or in the evaluation of special
functions. Unfortunately, the resulting sequences or series often do not
converge fast enough to be practically useful, or they even diverge,
which means that summation techniques have to be employed to obtain
numerically useful information. In such situations, it is often helpful
to employ so-called sequence transformations.

Formally, a sequence transformation $\mathcal{T}$ is simply a map
\begin{equation}
  \mathcal{T} \colon \{ S_{n} \} \to \{ T_{n} \} \, ,
   \qquad \ n \in \mathbb{N}_{0} \, ,
\end{equation}
which transforms a given sequence $\{ S_{n} \}$, whose convergence may be
unsatisfactory, to a new sequence $\{T_{n} \}$ with hopefully better
numerical properties. Of course, there is the minimal requirement that
the transformed sequence $\{T_{n} \}$ must have the same (generalized)
limit $S$ as the original sequence $\{ S_{n} \}$, but otherwise we have a
lot of freedom. It is common to say that $ \mathcal{T}$ accelerates
convergence if $\{T_{n} \}$ converges more rapidly to $S$ than $\{ S_{n}
\}$ according to
\begin{equation}
  \lim\limits_{n\to\infty} \, \frac{T_{n}-S}{S_{n}-S} \; = \; 0 \, .
\end{equation}

Many practically relevant sequences can be classified according to some
simple convergence types. Let us assume that the elements of a sequence
$\{ S_{n} \}$ with (generalized) limit $S$ satisfy the following
criterion resembling the well known ratio test:
\begin{equation}
  \label{LinLogConv}
  \lim_{n \to \infty} \, \frac {S_{n+1} - S} {S_{n} - S}
   \; = \; \rho \, .
\end{equation}
If $0 < \vert \rho \vert < 1$ holds, $\{ S_{n} \}$ converges
\emph{linearly}, if $\rho = 1$ holds, $\{ S_{n} \}$ converges
\emph{logarithmically}, and if $\rho = 0$ holds, it converges
\emph{hyperlinearly}. Obviously, $\vert \rho \vert > 1$ implies
divergence. The partial sums $\sum_{k=0}^n z^k = [1-z^{n+1}]/[1-z]$ of
the geometric series converge linearly for $\vert z \vert < 1$ as
$n \to \infty$, and they diverge for $\vert z \vert \ge 1$. The partial
sums $\sum_{k=0}^{n} (k+1)^{-s}$ of the Dirichlet series
$\zeta (s) = \sum_{\nu=0}^{\infty} (\nu+1)^{-s}$ for the Riemann zeta
function converge logarithmically for $\Re (s) > 1$ (see for example
\citep[Eq.\ (2.6)]{Weniger/Kirtman/2003}), and the partial sums
$\sum_{k=0}^{n} z^{k}/k!$ of the power series for $\exp (z)$ converge
hyperlinearly for $z \in \mathbb{C}$.

A sequence transformation $\mathcal{T}$ corresponds to an infinite set of
doubly indexed quantities $T_k^{(n)}$ with $k, n \in \mathbb{N}_{0}$.
Each $T_k^{(n)}$ is computed from a \emph{finite} subset
$\bigl\{ S_{n}, S_{n+1}, \dots, S_{n+\ell(k)} \bigr\}$ of
\emph{consecutive} elements, where $\ell = \ell (k)$. In our notation,
$T_{0}^{(n)}$ always corresponds to an untransformed sequence element
$T_{0}^{(n)} = S_{n}$.

The \emph{transformation order} $k$ is a measure for the complexity of
$T_k^{(n)}$, and the superscript $n$ corresponds to the minimal index $n$
of the string $\bigl\{ S_{n}, S_{n+1}, \dots, S_{n+\ell(k)} \bigr\}$ of
input data.  An increasing value of $k$ implies that the complexity of
the transformation process as well as $\ell = \ell (k)$ increases. Thus,
the application of a sequence transformation $\mathcal{T}$ to
$\{ S_{n} \}$ produces for every $k, n \ge \mathbb{N}_{0}$ a new
transformed element
\begin{equation}
  \label{Def:Tkn}
  T_k^{(n)} \; = \;
   T_k^{(n)} \bigl( S_{n}, S_{n+1}, \ldots , S_{n+\ell(k)} \bigr) \, .
\end{equation}

Sequence transformations are at least as old as calculus, and in
rudimentary form they are even much older. For the older history of
sequence transformations, we recommend a monograph by
\citet{Brezinski/1991a}, which discusses earlier work on continued
fractions, {P}ad\'{e} approximants, and related topics starting from the
17th century until 1945, as well as an article by \citet{Brezinski/2009},
which emphasizes pioneers of extrapolation methods. More recent
developments are treated in articles by
\citet{Brezinski/1996,Brezinski/2000c,Brezinski/2019} and
\citet{Brezinski/RedivoZaglia/2019}. A very complete list of older
references up to 1991 can be found in a bibliography by
\citet{Brezinski/1991b}.

Ever since the pioneering work of \citet{Shanks/1949} and
\citet{Wynn/1956a}, respectively, \emph{nonlinear} and in general also
\emph{nonregular} sequence transformations have dominated both research
and practical applications. For more details, see the monographs by
\citet{Brezinski/1977,Brezinski/1978,Brezinski/1980a},
\citet*{Brezinski/RedivoZaglia/1991a}, \citet*{Cuyt/Wuytack/1987},
\citet{Delahaye/1988}, \citet*{Liem/Lu/Shih/1995},
\citet*{Marchuk/Shaidurov/1983}, \citet{Sidi/2003}, \citet{Walz/1996},
and \citet{Wimp/1981}, or articles by
\citet*{Caliceti/Meyer-Hermann/Ribeca/Surzhykov/Jentschura/2007},
\citet{Homeier/2000a}, and \citet{Weniger/1989}. Sequence transformations
are also treated in the book by \citet*{Baker/Graves-Morris/1996} on
Pad\'{e} approximants, which can be viewed to be just a special class of
sequence transformations since they convert the partial sums of a power
series to a doubly indexed sequence of rational functions.

Practical applications of sequence transformations are emphasized in a
book by \citet*[Appendix A]{Bornemann/Laurie/Wagon/Waldvogel/2004} on
extreme digit hunting, in the most recent (third) edition of the book
\emph{Numerical Recipes} \citep{Press/Teukolsky/Vetterling/Flannery/2007}
(a review of the treatment of sequence transformations in \emph{Numerical
  Recipes} by \citet{Weniger/2007d} can be downloaded from the Numerical
Recipes web page), in a book by \citet*{Gil/Segura/Temme/2007} on the
evaluation of special functions (compare also the related review articles
by \citet*{Gil/Segura/Temme/2011} and by \citet{Temme/2007}), in the
recently published NIST Handbook of Mathematical Functions \citep[Chapter
3.9 Acceleration of Convergence]{Olver/Lozier/Boisvert/Clark/2010}, and
in a recent book by \citet{Trefethen/2013} on approximation
theory. Readers interested in current research both on and with sequence
transformations might look into the proceedings of a recent conference
\citep{Brezinski/RedivoZaglia/Weniger/2010a,Weniger/2010a}.

Certain sequence transformations are closely related to dynamical
systems. For example, Wynn's epsilon \citep{Wynn/1956a} and rho
\citep{Wynn/1956b} algorithm can be viewed to be just fully discrete
integrable systems
\citep{Nagai/Satsuma/1995,Papageorgiou/Grammaticos/Ramani/1993}. New
sequence transformations were also derived via this connection with
dynamical systems
\citep{Brezinski/He/Hu/Redivo-Zaglia/Sun/2012,He/Hu/Sun/Weniger/2011,Chang/He/Hu/Li/2018,
  Sun/Chang/He/Hu/2013}.

It is the basic assumption of all the commonly used sequence
transformations that the elements $S_{n}$ of a given sequence
$\{ S_{n} \}$ can be partitioned into a (generalized) limit $S$ and a
remainder $R_{n}$ according to
\begin{equation}
  S_{n} \; = \; S \, + \, R_{n} \, , \qquad n \in \mathbb{N}_{0} \, .
\end{equation}

How do sequence transformations accomplish an acceleration of convergence
or a summation of a divergent sequence? Let us assume that some
sufficiently good approximations $\{ \tilde{R}_{n} \}$ to the actual
remainders $\{ R_{n} \}$ of the elements of the sequence $\{ S_{n} \}$
are known. Elimination of these approximations $\tilde{R}_{n}$ yields a
new sequence
\begin{equation}
  S'_{n} \; = \; S_{n} - \tilde{R}_{n} \; = \; S + R_{n} -  \tilde{R}_{n}
   \; = \; S \, + \, R'_{n} \, .
\end{equation}
If $\tilde{R}_{n}$ is a good approximation to $R_{n}$, the transformed
remainder $R'_{n} = R_{n} - \tilde{R}_{n}$ vanishes more rapidly than the
original remainders $R_{n}$ as $n \to \infty$. At least conceptually,
this is what a sequence transformation tries to accomplish, although it
may actually use a completely different computational algorithm. Thus,
when trying to construct a new sequence transformation, we have to look
for arithmetic operations that lead in each step to an improved
approximate elimination of the truncation errors.

The idea of viewing a sequence transformation as an approximate
elimination procedure for the truncation errors becomes particularly
transparent in the case of model sequences whose remainders $R_{n}$
consists of a single exponential term:
\begin{equation}
  \label{AitModSeq}
  S_{n} \; = \; S \, + \, C \lambda^n \; , \qquad C \ne 0 \, , \,
   \quad \vert \lambda \vert < 1 \, , \, \quad n \in \mathbb{N}_{0} \, .
\end{equation}
This is an almost trivially simple model problem. Nevertheless, many
practically relevant sequences are known where exponential terms of the
type of $C \lambda^n$ form the dominant parts of their truncation errors.

A short computation shows that the so-called $\Delta^{2}$ formula (see for
example \citep[Eq.\ (5.1-4)]{Weniger/1989})
\begin{equation}
  \label{AitFor_1}
  \mathcal{A}_{1}^{(n)} \; = \; S_{n} \, - \,
   \frac{[\Delta S_{n}]^2}{\Delta^2 S_{n}} \, ,
    \qquad n \in \mathbb{N}_{0} \, ,
\end{equation}
is exact for the model sequence \eqref{AitModSeq}, i.e., we have
$\mathcal{A}_{1}^{(n)} = S$. Here, $\Delta$ is the finite difference
operator defined by $\Delta f (n) = f (n+1) - f (n)$.

Together with its numerous mathematically, but not necessarily
numerically equivalent explicit expressions (see for example \cite[Eqs.\
(5.1-6) - (5.1-12)]{Weniger/1989}), the $\Delta^{2}$ formula
\eqref{AitFor_1} is commonly attributed to \citet{Aitken/1926}, although
it is in fact much older. \citet[pp.\ 90 - 91]{Brezinski/1991a} mentioned
that in 1674 Seki Kowa, the probably most famous Japanese mathematician
of that period, tried to obtain better approximations to $\pi$ with the
help of the $\Delta^2$ formula (see also Osada's article \citep[Section
5]{Osada/2012}), and according to \citet[p.\ 5]{Todd/1962b}, the
$\Delta^{2}$ formula was in principle known to \citet{Kummer/1837}
already in 1837.

The power and practical usefulness of the $\Delta^{2}$ formula
\eqref{AitFor_1} is obviously limited since it only eliminates a single
exponential term. However, the $\Delta^{2}$ formula \eqref{AitFor_1} can
be iterated by using the output data $\mathcal{A}_{1}^{(n)}$ as input
data in the $\Delta^2$ formula. Repeating this process yields the
following nonlinear recursive scheme (see for example \cite[Eq.\
(5.1-15)]{Weniger/1989}):
\begin{subequations}
  \label{It_Aitken}
  \begin{align}
    \label{It_Aitken_a}
    \mathcal{A}_{0}^{(n)} & \; = \; S_{n} \, ,
     \qquad n \in \mathbb{N}_{0} \, ,
    \\
    \label{It_Aitken_b}
    \mathcal{A}_{k+1}^{(n)} & \; = \; \mathcal{A}_{k}^{(n)} - \frac
     {\bigl[\Delta \mathcal{A}_{k}^{(n)}\bigr]^2}
      {\Delta^2 \mathcal{A}_{k}^{(n)}} \, ,
       \qquad k, n \in \mathbb{N}_{0} \, .
  \end{align}
\end{subequations}
Here, $\Delta$ only acts on the superscript $n$ but not on the subscript
$k$ according to
$\Delta \mathcal{A}_{k}^{(n)} = \mathcal{A}_{k}^{(n+1)} -
\mathcal{A}_{k}^{(n)}$
(if necessary, we write $\Delta_{n} \mathcal{A}_{k}^{(n)}$ for the sake
of clarity).

The iterated $\Delta^{2}$ process \eqref{It_Aitken} is a fairly powerful
sequence transformation (see for example \cite[Table 13-1 on p.\ 328 or
Section 15.2]{Weniger/1989}). In \citep[p.\ 228]{Weniger/1989}, one finds
the listing of a FORTRAN 77 program that computes the iterated
$\Delta^{2}$ process using a single one-dimensional array.

There is an extensive literature on both the $\Delta^{2}$ formula
\eqref{AitFor_1} and its iteration \eqref{It_Aitken} (see for example
\citep{Brezinski/RedivoZaglia/2019,Weniger/2000b} and references
therein). But it is even more important that Aitken paved the way for the
Shanks transformation \citep{Shanks/1955} and for Wynn's epsilon
algorithm \citep{Wynn/1956a}, which permits a convenient recursive
computation of the Shanks transformation and thus also of Pad\'{e}
approximants. In 1956, Wynn also derived his so-called rho algorithm
\citep{Wynn/1956b}, whose recursive scheme closely resembles that of the
epsilon algorithm. These well known facts will be reviewed in more detail
in \cref{Sec:WynnEpsRhoAlg}.

The aim of our article is the clarification of the relationships among
known generalizations of Wynn's epsilon \citep{Wynn/1956a} and rho
\citep{Wynn/1956b} algorithms, and the construction of new
generalizations. We accomplish this by considering finite difference
equations in the transformation order $k$. To the best of our knowledge,
this is a novel approach. It leads to a better understanding of the
epsilon and rho algorithm and of many of their generalizations. In
addition, our approach also produces several new sequence
transformations.

It would be overly optimistic to assume that this article, which is the
first one to describe and apply our novel approach, could provide
definite answers to all questions occurring in this context. We also do
not think that the new transformations, which are presented in this
article, exhaust the potential of our novel approach.

The details of our novel approach will be discussed in
\cref{Sec:DiffEqTransOrd}. In \cref{Sec:NumEx}, we will show the
performance of our algorithms by some numerical examples. Finally, \cref{Sec:ConclDiscuss} is
devoted to discussions and conclusions.

%
\typeout{==> Section: Review of Wynn's epsilon and rho algorithms and of
  related transformations}
\section{Review of Wynn's epsilon and rho algorithm and of related
  transformations}
\label{Sec:WynnEpsRhoAlg}
%

%
\typeout{==> Subsection: The Shanks transformation and Wynn's epsilon algorithm}
\subsection{The Shanks transformation and Wynn's epsilon algorithm}
\label{Sub:ShanksTrWynnsEpsAlg}
%

An obvious generalization of the model sequence \eqref{AitModSeq}, which
leads to the $\Delta^{2}$ formula \eqref{AitFor_1}, is the following one,
which contains $k$ exponential terms instead of one:
\begin{equation}
  \label{ModSeqShanks_3}
  S_{n} \; = \; S + \sum_{j=0}^{k-1} \, C_{j} (\lambda_{j})^{n} \, ,
   \qquad k, n \in \mathbb{N}_{0} \, .
\end{equation}
It is usually assumed that the $\lambda_{j}$ are distinct
($\lambda_{m} \ne \lambda_{n}$ for $m \ne n$), and that they ordered
according to magnitude
($\vert \lambda_{0} \vert > \vert \lambda_{1} \vert > \dots > \vert
\lambda_{k-1} \vert$).

Although the $\Delta^2$ formula \eqref{AitFor_1} is by construction exact
for the model sequence \eqref{AitModSeq}, its iteration \eqref{It_Aitken}
with $k > 1$ is not exact for the model sequence
\eqref{ModSeqShanks_3}. Instead, the Shanks transformation, which is
defined by the following ratio of Hankel determinants \citep[Eq.\
(2)]{Shanks/1955},
\begin{equation}
  \label{ShanksTrDetRep}
  e_{k} (S_{n}) \; = \; \frac
   {\begin{vmatrix}
    S_{n} & S_{n+1} & \ldots & S_{n+k} \\
    \Delta S_{n} & \Delta S_{n+1} & \dots & \Delta S_{n+k} \\
     \vdots & \vdots & \ddots & \vdots \\
      \Delta S_{n+k-1} &\Delta S_{n+k} & \ldots &
       \Delta S_{n + 2 k - 1}
   \end{vmatrix}}
   {\begin{vmatrix}
    1 & 1 & \ldots & 1 \\
     \Delta S_{n} & \Delta S_{n+1} & \ldots & \Delta S_{n+k} \\
      \vdots & \vdots & \ddots & \vdots \\
       \Delta S_{n+k-1} & \Delta S_{n+k} & \ldots &
        \Delta S_{n + 2 k - 1}
   \end{vmatrix}} \, ,
    \qquad k, n \in \mathbb{N}_{0} \, .
\end{equation}
is exact for the model sequence \eqref{ModSeqShanks_3}. Actually, the
Shanks transformation had originally been derived in
\citeyear{Schmidt/1941} by \citet{Schmidt/1941}, but
\citeauthor{Schmidt/1941}'s article had largely been ignored.

Shanks' discovery of his transformation \eqref{ShanksTrDetRep} had an
enormous impact, in particular since he could show that it produces
Pad\'{e} approximants. Let
\begin{equation}
  \label{ParSum_PowSer}
  f_{n} (z) \; = \; \sum_{\nu=0}^{n} \, \gamma_{\nu} \, z^{\nu} \, ,
   \qquad n \in \mathbb{N}_{0} \, ,
\end{equation}
stand for the partial sums of the (formal) power series for some function
$f$. Then \cite[Theorem VI on p.\ 22]{Shanks/1955},
\begin{equation}
  \label{Shanks2Pade}
  e_k \bigl( f_{n} (z) \bigr) \; = \; [ k + n / k ]_f \, (z) \, ,
   \qquad k, n \in \mathbb{N}_{0} \, .
\end{equation}
We use the notation of \citet*[Eq.\ (1.2)]{Baker/Graves-Morris/1996} for
Pad\'{e} approximants.

The Hankel determinants in Shanks' transformation \eqref{ShanksTrDetRep}
can be computed recursively via the following non-linear recursion (see
for example \citep[p.\ 80]{Brezinski/RedivoZaglia/1991a}:
\begin{subequations}
  \label{HankelDetRec}
  \begin{align}
    \label{HankelDetRec_a}
      & H_{0} (u_{n}) \; = \; 1 \, , 
       \qquad H_{1} (u_{n}) \; = \; u_{n} \, ,
        \qquad n \in \mathbb{N}_{0} \, ,
    \\
    \label{HankelDetRec_b}
     & H_{k+2} (u_{n}) H_{k} (u_{n+2}) \; = \;  
      H_{k+1} (u_{n}) H_{k+1} (u_{n+2}) \, - \,
       \left[ H_{k+1} (u_{n+1}) \right]^{2} \, , 
        \qquad k, n \in \mathbb{N}_{0} \, .
  \end{align}
\end{subequations}
This recursive scheme is comparatively complicated. Therefore,
\citeauthor{Wynn/1956a}'s discovery of his celebrated epsilon algorithm
was a substantial step forwards \citep[Theorem on p.\ 91]{Wynn/1956a}:
\begin{subequations}
  \label{eps_al}
  \begin{align}
    \label{eps_al_a}
    \varepsilon_{-1}^{(n)} & \; = \; 0 \, ,
     \qquad \varepsilon_{0}^{(n)} \, = \, S_{n} \, ,
      \qquad  n \in \mathbb{N}_{0} \, , \\
    \label{eps_al_b}
    \varepsilon_{k+1}^{(n)} & \; = \; \varepsilon_{k-1}^{(n+1)} \, + \,
     \frac{1}{\varepsilon_{k}^{(n+1)} - \varepsilon_{k}^{(n)}} \, ,
      \qquad k, n \in \mathbb{N}_{0} \, ,
  \end{align}
\end{subequations}
\citet[Theorem on p.\ 91]{Wynn/1956a} showed that the
epsilon elements with \emph{even} subscripts give Shanks' transformation
\begin{equation}
  \label{eps_even_Shanks}
  \varepsilon_{2 k}^{(n)} \; = \; e_{k} (S_{n}) \, ,
   \qquad k,n \in \mathbb{N}_{0} \, ,
\end{equation}
whereas the elements with \emph{odd} subscripts are only auxiliary
quantities
\begin{equation}
  \label{eps_odd_Shanks}
  \varepsilon_{2 k + 1}^{(n)} \, = \; 1 / e_{k} ( \Delta S_{n} ) \, ,
   \qquad k,n \in \mathbb{N}_{0} \, ,
\end{equation}
which diverge if the transforms \eqref{eps_even_Shanks} converge to the
(generalized) limit $S$ of $\{ S_{n} \}$.

The epsilon algorithm is a close relative of Aitken's iterated
$\Delta^{2}$ process. We have
$\mathcal{A}_{1}^{(n)} = \varepsilon_{2}^{(n)}$, but for $k > 1$ we in
general have $\mathcal{A}_{k}^{(n)} \ne \varepsilon_{2k}^{(n)}$.
Nevertheless, the iterated $\Delta^{2}$ process and the epsilon algorithm
often have similar properties in convergence acceleration and summation
processes.

Because of the rhombus structure of the four-term recursion
\eqref{eps_al_b}, it appears that a program for Wynn's epsilon algorithm
would need either a single two-dimensional or at least two
one-dimensional arrays. However, \citet{Wynn/1965} could show with the
help of his \emph{moving lozenge technique}, which is described in a very
detailed way in \citep[Chapter 4.2.1.2]{Brezinski/1978}, that a single
one-dimensional array plus two auxiliary variables suffice. In
\citep[Chapter 4.3.2]{Brezinski/1978}, one finds listings of several
FORTRAN 66 programs for the epsilon algorithm.

As shown in \citep[Section 4.3]{Weniger/1989}, a slight improvement of
Wynn's moving lozenge technique is possible. The listing of a compact
FORTRAN 77 program for the epsilon algorithm using this modification can
be found in \citep[p.\ 222]{Weniger/1989}. In the book by \citet*[p.\
213]{Press/Teukolsky/Vetterling/Flannery/2007}, one finds a translation
of this FORTRAN 77 program to \texttt{C}.

It follows from \cref{Shanks2Pade,eps_even_Shanks} that the epsilon
algorithm produces Pad\'{e} approximants if its input data $S_{n}$ are
the partial sums \eqref{ParSum_PowSer} of a power series:
\begin{equation}
  \label{Eps_Pade}
  \varepsilon_{2 k}^{(n)} \; = \; [ k+n / k ]_f (z) \, ,
   \qquad k, n \in \mathbb{N}_{0} \, .
\end{equation}
Consequently, the epsilon algorithm is treated in books on Pad\'e
approximants such as the one by \citet{Baker/Graves-Morris/1996}, but
there is also an extensive literature dealing directly with it. According
to \citet[p.\ 120]{Wimp/1981}, over 50 articles on the epsilon algorithm
were published by Wynn alone, and at least 30 articles by Brezinski. As a
fairly complete source of references, Wimp recommends Brezinski's first
book \cite{Brezinski/1977}. However, this book was published in 1977, and
since then many more books and articles dealing with the epsilon
algorithm have been published. Any attempt of providing a reasonably
complete bibliography would be beyond the scope of this article.

The epsilon algorithm is not limited to input data that are the partial
sums of a (formal) power series. Therefore, it is actually more general
than Pad\'{e} approximants. For example, sequences of vectors can be used
as input data. For a relatively recent review, see
\citep{Graves-Morris/Roberts/Salam/2000a}. Both Shanks' transformation
and the epsilon algorithm are discussed in the NIST Handbook of
Mathematical Functions \citep*[Chapter 3.9(iv) Shanks'
Transformation]{Olver/Lozier/Boisvert/Clark/2010}.

Several generalizations of the epsilon algorithm are known. In 1972,
\citet[pp.\ 72 and 78]{Brezinski/1972} introduced what he called his
first and second generalization of the epsilon algorithm (see also
\citep[Chapter 2.6]{Brezinski/RedivoZaglia/1991a}). \citet[Algorithm A1
on p.\ 180]{Sablonniere/1991} and \citet[Algorithm 1 and Theorem 1 on p.\
255]{Sedogbo/1990} introduced other generalizations of the epsilon
algorithms that were specially designed to be effective for certain
logarithmically convergent interation sequences. Other generalizations
were considered by \citet*[Eq.\ (2.1)]{VandenBroeck/Schwartz/1979} and by
\citet*[Eqs.\ (43) - (44)]{Barber/Hamer/1982}.

%
\typeout{==> Subsection: Wynn's rho algorithm}
\subsection{Wynn's rho algorithm}
\label{Sub:WynnRhoAlg}
%

The iterated $\Delta^2$ process \eqref{It_Aitken} and the epsilon
algorithm \eqref{eps_al} are powerful accelerators for linearly
convergent sequences, and they are also able to sum many alternating
divergent series (see for example \citep[Section 15.2]{Weniger/1989} or
the rigorous convergence analysis of the summation of the Euler series in
\citep{Borghi/Weniger/2015}). However, both fail in the case of
logarithmic convergence \citep[Theorem 12]{Wynn/1966a}. Fortunately, in
\citeyear{Wynn/1956b} \citeauthor{Wynn/1956b} also derived his rho
algorithm \citep[Eq.\ (8)]{Wynn/1956b}
\begin{subequations}
\label{RhoAl}
\begin{align}
  \label{RhoAl_a}
  \rho_{-1}^{(n)} & \; = \; 0 \, ,
   \qquad \rho_{0}^{(n)} \; = \; S_{n} \, ,
    \qquad n \in \mathbb{N}_{0} \, ,
  \\
  \label{RhoAl_b}
  \rho_{k+1}^{(n)} & \; = \; \rho_{k-1}^{(n+1)} \, + \, \frac { x_{n+k+1}
    - x_{n}} {\rho_{k}^{(n+1)} - \rho_{k}^{(n)}} \, ,
     \qquad k, n \in \mathbb{N}_{0} \, ,
\end{align}
\end{subequations}
which is -- as for example emphasized by \citet{Osada/1996a} -- an
effective accelerator for many logarithmically convergent sequences.

As in the case of the formally almost identical epsilon algorithm, only
the elements $\rho_{2k}^{(n)}$ with \emph{even} subscripts provide
approximations to the limit of the input sequence. The elements
$\rho_{2k+1}^{(n)}$ with \emph{odd} subscripts are only auxiliary
quantities which diverge if the transformation process
converges. Actually, the elements $\rho_{2k}^{(n)}$ correspond to the
terminants of an interpolating continued fraction with interpolation
points $\{ x_{n} \}$ that are extrapolated to infinity (see for example
\cite[Chapter IV.1.4]{Cuyt/Wuytack/1987}). Consequently, the
interpolation points $\{ x_{n} \}$ must be positive, strictly increasing,
and unbounded:
\begin{subequations}
  \label{InterPolRgoAl}
  \begin{gather}
    \label{InterPolRgoAl_a}
    0 < x_{0} < x_{1} < \cdots < x_{m} < x_{m+1} < \cdots \, ,
    \\
    \label{InterPolRgoAl_b}
    \lim_{n \to \infty} x_{n} \; = \; \infty \, .
  \end{gather}
\end{subequations}

In the vast majority of all applications, Wynn's rho algorithm
\eqref{RhoAl} is used in combination with the interpolation points
$x_{n} = n+1$, yielding its standard form (see for example \citep[Eq.\
(6.2-4)]{Weniger/1989}):
\begin{subequations}
\label{RhoAlStand}
\begin{align}
  \label{RhoAlStand_a}
  \rho_{-1}^{(n)} & \; = \; 0 \, ,
   \qquad \rho_{0}^{(n)} \; = \; S_{n} \, ,
    \qquad n \in \mathbb{N}_{0} \, ,
  \\
  \label{RhoAlStand_b}
  \rho_{k+1}^{(n)} & \; = \; \rho_{k-1}^{(n+1)} \, + \, \frac {k+1}
   {\rho_{k}^{(n+1)} - \rho_{k}^{(n)}} \, ,
    \qquad k, n \in \mathbb{N}_{0} \, .
\end{align}
\end{subequations}

The standard form \eqref{RhoAlStand} is not able to accelerate the
convergence of \emph{all} logarithmically convergent sequences of
interest. This can be demonstrated by considering the following class of
logarithmically convergent model sequences:
\begin{equation}
  \label{ModSeqAlpha}
  S_{n} \; = \; S \, + \, (n+\beta)^{-\theta} \,
   \sum_{j=0}^{\infty} c_j / (n+\beta)^j \, ,
    \qquad n \in \mathbb{N}_{0} \, .
\end{equation}
Here, $\theta$ is a positive decay parameter and $\beta$ is a positive
shift parameter. The elements of numerous practically relevant
logarithmically convergent sequences $\{ S_{n} \}$ can at least in the
asymptotic domain of large indices $n$ be represented by series
expansions of that kind.

\citet[Theorem 3.2]{Osada/1990a} showed that the standard form
\eqref{RhoAlStand} accelerates the convergence of sequences of the type
of \cref{ModSeqAlpha} if the decay parameter $\theta$ is a positive
integer, but it fails if $\theta$ is non-integral.

If the decay parameter $\theta$ of a sequence of the type of
\cref{ModSeqAlpha} is explicitly known, \citeauthor{Osada/1990a}'s
variant of the rho algorithm can be used \citep[Eq.\ (3.1)]{Osada/1990a}:
\begin{subequations}
\label{OsRhoAl}
\begin{align}
  \label{OsRhoAl_a}
  \bar{\rho}_{-1}^{(n)} & \; = \; 0 \, , \qquad
   \bar{\rho}_{0}^{(n)} \; = \; S_{n} \, ,
    \qquad n \in \mathbb{N}_{0} \, ,
  \\
  \label{OsRhoAl_b}
  \bar{\rho}_{k+1}^{(n)} & \; = \; \bar{\rho}_{k-1}^{(n+1)} \, + \,
   \frac
    {k+\theta} {\bar{\rho}_{k}^{(n+1)} - \bar{\rho}_{k}^{(n)}} \, ,
     \qquad k, n \in \mathbb{N}_{0} \, .
\end{align}
\end{subequations}
\citeauthor{Osada/1990a} showed that his variant \eqref{OsRhoAl}
accelerates the convergence of sequences of the type of
\cref{ModSeqAlpha} for known, but otherwise arbitrary, i.e., also for
non-integral $\theta > 0$. He proved the following asymptotic estimate
for fixed $k$ \citep[Theorem 4]{Osada/1990a}:
\begin{equation}
  {\bar \rho}_{2 k}^{(n)} \, - \, S \; = \;
   \mathrm{O} \bigl( n^{-\theta - 2k} \bigr) \, ,
    \qquad n \to \infty \, .
\end{equation}

As remarked above, the decay parameter $\theta$ must be explicitly known
if \citeauthor{Osada/1990a}'s transformation \eqref{OsRhoAl} is to be applied to a sequence
of the type of \cref{ModSeqAlpha}. An approximation to $\theta$ can be
obtained with the help of the following nonlinear transformation, which
was first derived in a somewhat disguised form by \citet[p.\
419]{Drummond/1976} and later rederived by \citet*[Eq.\
(4.1)]{Bjoerstad/Dahlquist/Grosse/1981}:
\begin{equation}
  \label{DecPar}
  T_{n} \; = \; \frac
   { [\Delta^2 S_{n} ] \, [\Delta^2 S_{n+1} ]}
    {[\Delta S_{n+1} ] \, [\Delta^2 S_{n+1} ] \, - \,
     [\Delta S_{n+2} ] \, [\Delta^2 S_{n}]} \, - \, 1 \, ,
      \qquad n \in \mathbb{N}_{0} \, .
\end{equation}
\citet*[Eq.\ (4.1)]{Bjoerstad/Dahlquist/Grosse/1981} showed that
\begin{equation}
  \label{T_n}
  \theta \; = \; T_{n} \, + \, O (1/n^2) \,,\qquad n \to \infty \, ,
\end{equation}
if the input data are the elements of a sequence of the type of
\cref{ModSeqAlpha}.

The variants of Wynn's rho algorithm can also be iterated in the spirit
of Aitken's iterated $\Delta^{2}$ process \eqref{It_Aitken}, which can
also be viewed to be an iteration of $\varepsilon_{2}^{(n)}$. By
iterating the expression for $\rho_{2}^{(n)}$ involving unspecified
interpolation points $x_{n}$, a rho analog of Aitken's iterated
$\Delta^{2}$ was constructed \cite[Eq.\ (2.10)]{Weniger/1991}. An
alternative iteration was derived by \citet*[Eq.\
(2.25)]{Bhowmick/Bhattacharya/Roy/1989}, which is, however, significantly
less efficient than the iteration \cite[Eq.\ (2,10)]{Weniger/1991}
involving unspecified interpolation points (compare \cite[Table
I]{Weniger/1991}). This is caused by the fact that \citet*[Eq.\
(2.25)]{Bhowmick/Bhattacharya/Roy/1989} had started from the standard
form \eqref{RhoAlStand} and not from the general form
\eqref{RhoAl}. \citeauthor{Osada/1990a}'s variant \eqref{OsRhoAl} of the
rho algorithm can also be iterated. It yields a recursive scheme
\citep[Eq.\ (2.29)]{Weniger/1991}, which was originally derived by
\citet*[Eq.\ (2.4)]{Bjoerstad/Dahlquist/Grosse/1981} who called the
resulting algorithm a modified $\Delta^{2}$ formula.

Wynn's epsilon algorithm \eqref{eps_al} is undoubtedly the currently most
important and most often employed sequence transformation. Wynn's rho
algorithm \eqref{RhoAl} and its standard form \eqref{RhoAlStand},
\citeauthor{Osada/1990a}'s generalized rho algorithm \eqref{OsRhoAl}, as
well as the iterations mentioned above are also very important and
practically useful sequence transformations. Accordingly, there are
numerous references describing successful applications (far too many to
be cited here). In addition, these transformations were very
consequential from a purely theoretical point of view. They inspired the
derivation of numerous other sequence transformations and ultimately also
this article.

%
\typeout{==> Subsection: Brezinski's theta algorithm and its iteration}
\subsection{Brezinski's theta algorithm and its iteration}
\label{Sub:BrezThetaAlg}
%

Soon after their derivation it had become clear that epsilon and rho have
largely complementary properties. The epsilon algorithm is a powerful
accelerator of linear convergence and also frequently works well in the
case of divergent series, but it does not accomplish anything substantial
in the case of logarithmic convergence. In contrast, the rho algorithm is
ineffective in the case of linear convergence and even more so in the
case of divergent series, but is able to speed up the convergence of many
logarithmically convergent sequences.

This observation raised the question whether it would be possible to
modify the recursions of epsilon and/or rho in such a way that the
resulting transformation would possess the advantageous features of both
epsilon and rho. This goal was achieved in 1971 by Brezinski who suitably
modified the epsilon recursion \eqref{eps_al} to obtain his theta
algorithm \citep[p.\ 729]{Brezinski/1971a} (compare also \citep[Chapter
2.3]{Brezinski/RedivoZaglia/1991a} or \citep[Section
10.1]{Weniger/1989}):
\begin{subequations}
  \label{ThetaAl}
  \begin{align}
    \label{ThetaAl_a}
    \vartheta_{-1}^{(n)} & \; = \; 0 \, , \qquad
     \vartheta_{0}^{(n)} \; = \; S_{n} \, ,
      \qquad n \in \mathbb{N}_{0} \, ,
    \\
    \label{ThetaAl_b}
    \vartheta_{2 k + 1}^{(n)} & \; = \; \vartheta_{2 k-1}^{(n+1)}
     + \frac{1}{\Delta \vartheta_{2 k}^{(n)}} \, ,
      \qquad k, n \in \mathbb{N}_{0} \, ,
    \\
    \label{ThetaAl_c}
    \vartheta_{2 k+2}^{(n)} & \; = \; \vartheta_{2 k}^{(n+1)} + \frac
     {\bigl[ \Delta \vartheta_{2 k}^{(n+1)} \bigr] \,
      \bigl[\Delta \vartheta_{2 k + 1}^{(n+1)} \bigr]}
     {\Delta^{2} \vartheta_{2 k+1}^{(n)}} \, ,
      \qquad k, n \in \mathbb{N}_{0} \, .
  \end{align}
\end{subequations}
We again assume that the finite difference operator $\Delta$ only acts on
$n$ but not on $k$.

As in the case of Wynn's epsilon and rho algorithm, only the elements
$\vartheta_{2k}^{(n)}$ with even subscripts provide approximations to the
(generalized) limit of the sequence $\{ S_{n} \}$ to be transformed. The
elements $\vartheta_{2k+1}^{(n)}$ with odd subscripts are only auxiliary
quantities which diverge if the whole process converges.

Brezinski's derivation of his theta algorithm was purely experimental,
but it was certainly a very successful experiment. Extensive numerical
studies performed by \citet{Smith/Ford/1979,Smith/Ford/1982} showed that
among virtually all sequence transformations known at that time only the
theta algorithm together with Levin's $u$ and $v$ transformations
\citep[Eqs.\ (59) and (68)]{Levin/1973} provided consistently good
results for a very wide range of test problems. This positive verdict was
also confirmed by later numerical studies.

In \citeauthor{Brezinski/1978}'s second book \citep[pp.\ 368 -
370]{Brezinski/1978}, one finds the listing of a FORTRAN IV program which
computes \citeauthor{Brezinski/1978}'s theta algorithm \ref{ThetaAl}
using three one-dimensional arrays, and in \citep[pp.\ 279 -
281]{Weniger/1989}, one finds the listing of a FORTRAN 77 program which
computes theta using two one-dimensional arrays.

It is possible to proceed in the spirit of Aitken's iterated $\Delta^{2}$
process \eqref{It_Aitken} and to construct iterations of Brezinski's
theta algorithm \eqref{ThetaAl}. A natural starting point for such an
iteration would be $\vartheta_{2}^{(n)}$ which possesses many alternative
expressions (see for example \citep[Eqs.\ (10.3-1) -
(10.3-3)]{Weniger/1989} or \citep[Eq.\
(4.10)]{Weniger/2000b}):
\begin{align}
  \label{theta_2_1}
  \vartheta_{2}^{(n)} & \; = \; S_{n+1} \, - \, \frac
   {\bigl[ \Delta S_{n} \bigr] \bigl[ \Delta S_{n+1} \bigr]
    \bigl[ \Delta^{2} S_{n+1} \bigr]}
     {\bigl[ \Delta S_{n+2} \bigr] \bigl[ \Delta^{2} S_{n} \bigr] -
      \bigl[ \Delta S_{n} \bigr] \bigl[ \Delta^{2} S_{n+1} \bigr]}
  \\
  \label{theta_2_2}
  & \; = \; \frac
   {S_{n+1} \, [\Delta S_{n+2}] \, [\Delta^{2} S_{n}] -
    S_{n+2} \, [\Delta S_{n}] \, [\Delta^{2} S_{n+1}]}
     {[\Delta S_{n+2}] \, [\Delta^{2} S_{n}] -
      [\Delta S_{n}] \, [\Delta^{2} S_{n+1}]}
  \\
  \label{theta_2_3}
  & \; = \; \frac
  {\Delta^{2} \, [ S_{n+1} / \Delta S_{n} ]}
   {\Delta^{2} \, [1/\Delta S_{n}]}
  \\
  \label{theta_2_7}
  & \; = \; S_{n+3} - \frac{ \bigl[\Delta S_{n+2}
   \bigr] \Bigl\{ \bigl[ \Delta S_{n+2} \bigr]
    \bigl[ \Delta^{2} S_{n} \bigr] + \bigl[ \Delta S_{n+1} \bigr]^{2} -
     \bigl[ \Delta S_{n+2} \bigr] \bigl[ \Delta S_{n} \bigr] \Bigr\} }
      {\bigl[ \Delta S_{n+2} \bigr] \bigl[\Delta^{2} S_{n} \bigr] -
       \bigl[\Delta S_{n} \bigr] \bigl[\Delta^{2} S_{n+1} \bigr]} \, .
\end{align}
As in the case of $\mathcal{A}_{2}^{(n)}$, these expressions are
obviously mathematically equivalent, but they may differ -- possibly
substantially -- with respect to their numerical stability.

By iterating the explicit expression \eqref{theta_2_1} for
$\vartheta_{2}^{(n)}$ in the spirit of Aitken's iterated $\Delta^{2}$
process \eqref{It_Aitken}, the following nonlinear recursive system was derived \citep[Eq.\ (10.3-6)]{Weniger/1989}:
\begin{subequations}
  \label{ThetIt}
  \begin{align}
    \label{ThetIt_a}
    \mathcal{J}_{0}^{(n)} \; = \; & S_{n} \, ,
     \qquad n \in \mathbb{N}_{0} \, ,
    \\
    \label{ThetIt_b}
    \mathcal{J}_{k+1}^{(n)} \; = \; & \mathcal{J}_{k}^{(n+1)} - \frac
     {\bigl[ \Delta \mathcal{J}_{k}^{(n)} \bigr] \bigl[
      \Delta \mathcal{J}_{k}^{(n+1)} \bigr] \bigl[
       \Delta^2 \mathcal{J}_{k}^{(n+1)} \bigr]}
        {\bigl[ \Delta \mathcal{J}_{k}^{(n+2)} \bigr]
         \bigl[ \Delta^2 \mathcal{J}_{k}^{(n)} \bigr] -
          \bigl[ \Delta \mathcal{J}_{k}^{(n)} \bigr]
           \bigl[\Delta^2 \mathcal{J}_{k}^{(n+1)} \bigr]} \, ,
     \notag \\
     & \qquad k, n \in \mathbb{N}_{0} \, .
  \end{align}
\end{subequations}
In \citep[pp.\ 284 - 285]{Weniger/1989}, one finds the listing a FORTRAN
77 program that computes $\mathcal{J}_k^{(n)}$ using a single
one-dimensional array.

Based on the explicit expression \eqref{theta_2_7} for
$\vartheta_{2}^{(n)}$, an alternative recursive scheme for
$\mathcal{J}_k^{(n)}$ was derived in \citep[Eq.\ (4.11)]{Weniger/2000b}
and used for the prediction of unknown power series coefficients.

The iterated transformation $\mathcal{J}_k^{(n)}$ has similar properties
as the theta algorithm from which it was derived. Both are very powerful
as well as very versatile. $\mathcal{J}_k^{(n)}$ is not only an effective
accelerator for linear convergence as well as able to sum divergent
series, but it is also able to accelerate the convergence of many
logarithmically convergent sequences and series (see for example
\citep{Sidi/2002,Weniger/1989,Weniger/1991,Weniger/1994b,Weniger/2000b}
and references therein).

As discussed in \cref{App:LevTr}, the special case
$\vartheta_{2}^{(n)} = \mathcal{J}_{1}^{(n)}$ of the iterated
transformation $\mathcal{J}_k^{(n)}$ can be related to the $u$ and $v$
variants of Levin's general transformation \citep[Eq.\ (22)]{Levin/1973}
defined by \cref{GenLevTr}. For Levin's $u$ and $v$ variants \citep[Eqs.\
(59) and (68)]{Levin/1973}, we use the notation \citep[Eqs.\ (7.3-5) and
(7.3-11)]{Weniger/1989}:
\begin{align}
  \label{uLevTr}
  u_{k}^{(n)} (\beta, S_{n}) & \; = \;
   \mathcal{L}_{k}^{(n)} (\beta, S_{n}, (\beta + n) \Delta S_{n-1}) \, ,
  \\
  \label{vLevTr}
  v_{k}^{(n)} (\beta, s_{n}) & \; = \;
   \mathcal{L}_{k}^{(n)} (\beta, S_{n},
    \Delta S_{n-1} \Delta S_{n} / [\Delta S_{n-1} - \Delta S_{n}]) \, .
\end{align}
\Cref{theta_2_1,uLevTr_5,vLevTr_2} imply (see also \citep[Eq.\
(10.3-4)]{Weniger/1989}):
\begin{align}
  \label{theta_2<->u_2}
  \vartheta_{2}^{(n)} & \; = \; u_{2}^{(n+1)} (\beta, S_{n+1}) \, ,
  \\
  \label{theta_2<->v_1}
  & \; = \; v_{1}^{(n+1)} (\beta, S_{n+1}) \, ,
   \qquad n \in \mathbb{N}_{0} \, , \quad \beta > 0 \, .
\end{align}

%
\typeout{==> Section: Difference equations with respect to the
  transformation order}
\section{Difference equations with respect to the transformation order}
\label{Sec:DiffEqTransOrd}
%

%
\typeout{==> Subsection: General Considerations}
\subsection{General Considerations}
\label{Sub:GeneralConsiderations}
%

In this Section, we want to clarify the relationship between known
generalizations of Wynn's epsilon and rho algorithm, and -- more
important -- we also want to derive new generalizations. For that
purpose, we pursue a novel approach and look at their recursions from the
perspective of finite difference equations with respect to the
transformation order $k$.

Our approach is based on the observation that the recursions
\eqref{eps_al_b}, \eqref{RhoAl_b}, \eqref{RhoAlStand_b}, and
\eqref{OsRhoAl_b} for Wynn's epsilon and rho algorithm and its special
cases all possess the following general structure:
\begin{equation}
  \label{Tkn_2_Fkn}
  F_{k}^{(n)} \; = \; \left[ T_{k+1}^{(n)} - T_{k-1}^{(n+1)} \right] \,
   \left[ T_{k}^{(n+1)} - T_{k}^{(n)} \right] \, ,
    \qquad k, n \in \mathbb{N}_{0} \, .
\end{equation}
The quantities $F_{k}^{(n)}$, which in general depend on both $k$ and
$n$, fully characterize the corresponding transformations. The algebraic
properties of sequence transformations satisfying a recursion of the type
of \cref{Tkn_2_Fkn} had already been studied by Brezinski in his thesis
\citep[pp.\ 120 - 127]{Brezinski/1971b} (see also \citep[pp.\ 106 -
107]{Brezinski/RedivoZaglia/1991a}).

As shown later, many known variants of epsilon and rho can be classified
according to the properties of $F_{k}^{(n)}$. This alone makes our
subsequent analysis based on \cref{Tkn_2_Fkn} useful. But it is our
ultimate goal to construct new generalizations of epsilon and rho.

In principle, we could try to achieve our aim of constructing new
transformations by simply choosing new and hopefully effective
expressions for $F_{k}^{(n)}$ in \cref{Tkn_2_Fkn}. In combination with
suitable initial conditions, this would produce recursive schemes for new
sequence transformations. To the best of our knowledge, questions of that
kind have not been treated in the literature yet. So -- frankly speaking
-- we simply do not know how a practically useful $F_{k}^{(n)}$ should
look like.

In order to get some help in this matter, we will make some specific
assumptions about the dependence of the $F_{k}^{(n)}$ as functions of $k$
and $n$. These assumptions will allow us to investigate whether and under
which conditions the $F_{k}^{(n)}$ satisfy certain finite difference
equations in the transformation order $k$. It will become clear later
that our analysis of the resulting finite difference equations in $k$
leads to some insight which ultimately makes it possible to construct new
sequence transformations.

For an analysis of the \emph{known} variants of epsilon and rho, the full
generality of \cref{Tkn_2_Fkn} is only needed in the case of the general
form \eqref{RhoAl} of Wynn's rho algorithm with unspecified interpolation
points $\{ x_{n} \}$. In all other cases, the apparent $n$-dependence of
the right-hand side of \cref{Tkn_2_Fkn} cancels out and we only have to
deal with a simpler, $n$-independent expression of the following general
structure:
\begin{equation}
  \label{Tkn_2_Fk}
  F_{k} \; = \; \left[ T_{k+1}^{(n)} - T_{k-1}^{(n+1)} \right] \,
   \left[ T_{k}^{(n+1)} - T_{k}^{(n)} \right] \, ,
    \qquad k, n \in \mathbb{N}_{0} \, .
\end{equation}
In the case of the epsilon algorithm \eqref{eps_al}, we have $F_{k} = 1$,
in the case of the standard form \eqref{RhoAlStand} of the rho algorithm,
we have $F_{k} = k+1$, and in the case of Osada's generalization
\eqref{OsRhoAl}, we have $F_{k} = k+\theta$.

%
\typeout{==> Subsection: First order difference equations}
\subsection{First order difference equations}
\label{Sub:FirstOrderDifferenceEquation}
%

The simplest case $F_{k} = 1$ in \cref{Tkn_2_Fk}, which corresponds to
Wynn's epsilon algorithm \eqref{eps_al}, satisfies the following first
order finite difference equation in $k$:
\begin{equation}
  \label{Diff_Eq_1stOrd_k}
  \Delta_{k} F_{k} \; = \; F_{k+1} - F_{k} \; = \; 0 \, ,
   \qquad k \in \mathbb{N}_{0} \, .
\end{equation}
\Cref{Diff_Eq_1stOrd_k} implies that the recursion \eqref{Tkn_2_Fk} can
also be expressed as follows:
\begin{align}
  \label{Diff_Eq_1stOrd_k_Tkn}
  & \left[ T_{k+2}^{(n)} - T_{k}^{(n+1)} \right] \,
   \left[ T_{k+1}^{(n+1)} - T_{k+1}^{(n)} \right]
 \notag
 \\
 & \qquad - \,
  \left[ T_{k+1}^{(n)} - T_{k-1}^{(n+1)} \right] \,
   \left[ T_{k}^{(n+1)} - T_{k}^{(n)} \right] \; = \; 0 \, .
\end{align}
This reformulation of \cref{Tkn_2_Fk} does not look like a great
achievement. It also appears to be a bad idea to use the more complicated
recursion \eqref{Diff_Eq_1stOrd_k_Tkn} instead of the much simpler
recursion \eqref{Tkn_2_Fk}. However, the more complicated recursion
\eqref{Diff_Eq_1stOrd_k_Tkn} will help us to understand better certain
features of the simpler recursion \eqref{Tkn_2_Fk} and of the epsilon
algorithm \eqref{eps_al}. In addition, the recursion
\eqref{Diff_Eq_1stOrd_k_Tkn} as well as related expressions occurring
later open the way for the construction of new sequence transformations.

The first order finite difference equation $\Delta_{m} f (m) = f (m+1) -
f (m) = 0$ with $m \in \mathbb{N}_{0}$ possesses the general solution $f
(m) = f (0)$. Thus, \cref{Diff_Eq_1stOrd_k} implies $F_{k} = F_{0}$, and
\cref{Tkn_2_Fk} can be rewritten as follows:
\begin{equation}
  \label{Rec_Tkn_F0}
  T_{k+1}^{(n)} \; = \; T_{k-1}^{(n+1)} \, + \,
   \frac{F_{0}}{T_{k}^{(n+1)} - T_{k}^{(n)}} \, ,
    \qquad k, n \in \mathbb{N}_{0} \, .
\end{equation}
If we choose $F_{0} = 1$, we obtain Wynn's epsilon algorithm
\eqref{eps_al}.

Alternatively, we could also assume that $F_{0}$ is a non-zero
constant. By combining \cref{Rec_Tkn_F0} with the epsilon-type initial
conditions $T_{-1}^{(n)} = 0$ and $T_{0}^{(n)} = S_{n}$, we obtain
$T_{2k}^{(n)} = \varepsilon_{2k}^{(n)}$ and
$T_{2k+1}^{(n)} = F_{0} \, \varepsilon_{2k+1}^{(n)}$. But only the
transforms $T_{2k}^{(n)} = \varepsilon_{2k}^{(n)}$ with \emph{even}
subscripts provide approximations to the limit of the input sequence,
whereas the transforms $T_{2k+1}^{(n)} = F_{0} \varepsilon_{2k+1}^{(n)}$
with \emph{odd} subscripts are only auxiliary quantities. Thus, the
modified recursion \eqref{Rec_Tkn_F0} with $F_{0} \ne 0$ does not produce
anything new, and the epsilon algorithm remains invariant if we replace
in \cref{eps_al_b} $1/[\varepsilon_{k}^{(n+1)} - \varepsilon_{k}^{(n)}]$
by $F_{0}/[\varepsilon_{k}^{(n+1)} - \varepsilon_{k}^{(n)}]$. To the best
of our knowledge, this invariance of the epsilon algorithm had not been
observed before.

A similar invariance exists in the case of the general form \eqref{RhoAl}
of Wynn's rho algorithm and its variants. If we replace the interpolation
points $\{ x_{n} \}$ by the scaled interpolation points
$\{ \alpha x_{n} \}$ with $\alpha > 0$, we obtain
$\rho_{2k}^{(n)} (\alpha x_{n}) = \rho_{2k}^{(n)} (x_{n})$ and
$\rho_{2k+1}^{(n)} (\alpha x_{n}) = \alpha \rho_{2k+1}^{(n)} (x_{n})$.
This invariance can also be deduced from the determinantal representation
for $\rho_{2k}^{(n)}$ \citep[Theoreme 2.1.52 on p.\ 97]{Brezinski/1978}.

Our assumption, that $F_{k}$ satisfies the first order finite difference
equation \eqref{Diff_Eq_1stOrd_k} in $k$ and is therefore constant as a
function of $k$, is apparently too restrictive to allow the construction
of \emph{useful} alternative transformations beyond Wynn's epsilon
algorithm.

Let us now assume that the $n$-dependence of $F_{k}^{(n)}$ is such that
the first order finite difference equation
\begin{equation}
  \label{Diff_Eq_1stOrd_n-dep_k}
  \Delta_{k} F_{k}^{(n)} \; = \; F_{k+1}^{(n)} - F_{k}^{(n)}
   \; = \; 0 \, , \qquad k, n \in \mathbb{N}_{0} \, ,
\end{equation}
is satisfied for all $n \in \mathbb{N}_{0}$. Then, the finite difference
equation \eqref{Diff_Eq_1stOrd_n-dep_k} implies $F_{k}^{(n)} =
F_{0}^{(n)}$ and
\begin{equation}
  \label{Rec_Tkn_F0n}
  T_{k+1}^{(n)} \; = \; T_{k-1}^{(n+1)} \, + \,
   \frac {F_{0}^{(n)}} {T_{k}^{(n+1)} - T_{k}^{(n)}} \, ,
    \qquad k, n \in \mathbb{N}_{0} \, .
\end{equation}
If we choose $F_{0}^{(n)} = x_{n+1} - x_{n}$ in \cref{Rec_Tkn_F0n}, we
obtain Brezinski's first generalization of the epsilon algorithm
\citep[p.\ 72]{Brezinski/1972}, and if we choose $F_{0}^{(n)} = x_{n+k+1}
- x_{n+k}$, we formally obtain Brezinski's second generalization
\citep[p.\ 78]{Brezinski/1972}. However, Brezinski's second
generalization is only a special case of the recursion
\eqref{Rec_Tkn_F0n} if the original $k$-dependence of the interpolation
points $x_{n+k}$ is eliminated by forming the differences $\Delta x_{n+k}
= x_{n+k+1} - x_{n+k}$.

By combining the recursion \eqref{Rec_Tkn_F0n} with suitable initial
conditions and suitable choices for $F_{0}^{(n)}$, we obtain new sequence
transformations. Possible choices of the initial conditions will be be
discussed in \cref{Sub:ChoiceIniCon}.

%
\typeout{==> Subsection: Second order difference equations}
\subsection{Second order difference equations}
\label{Sub:SecondOrderDifferenceEquation}
%

We have $F_{k} = k+1$ in the case of the standard form \eqref{RhoAlStand}
of Wynn's rho algorithm, and in the case of Osada's generalization
\eqref{OsRhoAl} we have $F_{k} = k+\theta$. Thus, it looks like a natural
idea to assume that $F_{k}$ is linear in $k$, which means that $F_{k}$ is
annihilated by $\Delta_{k}^{2}$:
\begin{equation}
  \label{Diff_Eq_2ndOrd_k}
  \Delta_{k}^{2} F_{k} \; = \; F_{k+2} - 2 F_{k+1} + F_{k} \; = \; 0 \, ,
  \qquad k \in \mathbb{N}_{0} \, .
\end{equation}
The second order finite difference equation \eqref{Diff_Eq_2ndOrd_k} is
equivalent to the recursion
\begin{align}
  \label{Diff_Eq_2ndOrd_k_Tkn}
  & \left[ T_{k+3}^{(n)} - T_{k+1}^{(n+1)} \right] \,
   \left[ T_{k+2}^{(n+1)} - T_{k+2}^{(n)} \right]  \, - \,
    2 \, \left[ T_{k+2}^{(n)} - T_{k}^{(n+1)} \right] \,
     \left[ T_{k+1}^{(n+1)} - T_{k+1}^{(n)} \right]
  \notag
  \\
  & \qquad \, + \,
   \left[ T_{k+1}^{(n)} - T_{k-1}^{(n+1)} \right] \,
    \left[ T_{k}^{(n+1)} - T_{k}^{(n)} \right]
     \; = \; 0,  \qquad k, n \in \mathbb{N}_{0} \, ,
\end{align}
which is even more complicated than the analogous recursion
\eqref{Diff_Eq_1stOrd_k_Tkn} based on the first order difference equation
\eqref{Diff_Eq_1stOrd_k}.

The second order finite difference equation $\Delta_{m}^{2} f (m) = f
(m+2) - 2 f (m+1) + f (m) = 0$ with $m \in \mathbb{N}_{0}$ possesses the
general solution $f (m) = f (0) + m [f (1) - f (0)]$. Thus,
\cref{Diff_Eq_2ndOrd_k} implies $F_{k} = F_{0} + k [ F_{1} - F_{0}]$,
which means that \cref{Tkn_2_Fk} can be reformulated as follows:
\begin{equation}
  \label{Rec_Tkn_F0_F1}
  T_{k+1}^{(n)} \; = \; T_{k-1}^{(n+1)} \, + \,
   \frac{F_{0} + k [F_{1} - F_{0}]}{T_{k}^{(n+1)} - T_{k}^{(n)}} \, ,
    \qquad k, n \in \mathbb{N}_{0} \, .
\end{equation}
This recursion contains several variants of Wynn's rho algorithm as
special cases. If we choose $F_{0} = 1$ and $F_{1} = 2$, we obtain the
standard form \eqref{RhoAlStand} of Wynn's rho algorithm, and if we
choose $F_{0} = \theta$ and $F_{1} = \theta + 1$, we obtain Osada's
variant \eqref{OsRhoAl} of the rho algorithm. Finally,
$F_{0} = F_{1} \ne 0$ yields epsilon. Other choices for $F_{0}$ and
$F_{1}$ are obviously possible, but to the best of our knowledge they
have not been explored yet in the literature.

Our approach based on finite difference equations in the transformation
order $k$ was successful in the case of the variants of rho mentioned
above, but the general form \eqref{RhoAl} of Wynn's rho algorithm with
essentially arbitrary interpolation points $\{ x_{n} \}$ remains
elusive. In the case of essentially arbitrary interpolation points
$\{ x_{n} \}$, we cannot express the numerator $x_{n+k+1} - x_{n}$ of the
ratio in \cref{RhoAl_b} as an $n$-independent $F_{k}$ according to
\cref{Tkn_2_Fk}.

At least, we have to employ the more general expression \eqref{Tkn_2_Fkn}
involving a general $n$-dependent $F_{k}^{(n)}$. But even this does not
guarantee success. Of course, $F_{k}^{(n)} = x_{n+k+1} - x_{n}$ with
\emph{unspecified} and essentially \emph{arbitrary} interpolation points
$\{ x_{n} \}$ implies the finite difference equation
$\Delta_{k}^{m} F_{k}^{(n)} = \Delta_{k}^{m} x_{n+k+1}$ in $k$. But to
achieve any progress, we would have to make specific assumptions on the
index dependence of the interpolation points $\{ x_{n} \}$.

Therefore, let us now assume that the $n$-dependence of $F_{k}^{(n)}$ is
such that
\begin{equation}
  \label{Diff_Eq_2ndOrd_n-dep_k}
  \Delta_{k}^{2} F_{k}^{(n)} \; = \;
   F_{k+2}^{(n)} - 2 F_{k+1}^{(n)} + F_{k}^{(n)} \; = \; 0 \, ,
    \qquad k \in \mathbb{N}_{0} \, ,
\end{equation}
is satisfied for arbitrary $n \in \mathbb{N}_{0}$. This yields the
general solution $F_{k}^{(n)} = F_{0}^{(n)} + k [ F_{1}^{(n)} -
F_{0}^{(n)}]$, which implies that the recursion \eqref{Rec_Tkn_F0_F1} can
be generalized as follows:
\begin{equation}
  \label{Rec_Tkn_Fn0_Fn1}
  T_{k+1}^{(n)} \; = \; T_{k-1}^{(n+1)} \, + \, \frac
   {F_{0}^{(n)} + k [F_{1}^{(n)} - F_{0}^{(n)}]}
    {T_{k}^{(n+1)} - T_{k}^{(n)}} \, ,
     \qquad k, n \in \mathbb{N}_{0} \, .
\end{equation}
By combining this recursion with suitable initial conditions, we obtain a
whole class of new sequence transformations. Possible choices of the
initial conditions are discussed in \cref{Sub:ChoiceIniCon}.

%
\typeout{==> Subsection: Initial conditions for the new transformation}
\subsection{Initial conditions for the new transformation}
\label{Sub:ChoiceIniCon}
%

Let us first analyze the recursion \eqref{Diff_Eq_1stOrd_k_Tkn} which can
be reformulated as follows:
\begin{equation}
  \label{Diff_Eq_1stOrd_k_Tkn_1}
  T_{k+2}^{(n)} \; = \; T_{k}^{(n+1)} + \frac
   {\left[ T_{k+1}^{(n)} - T_{k-1}^{(n+1)} \right] \,
    \left[ T_{k}^{(n+1)} - T_{k}^{(n)} \right]}
     {T_{k+1}^{(n+1)} - T_{k+1}^{(n)}} \, .
\end{equation}
We have to investigate for which values of $k \in \mathbb{N}_{0}$ this
recursion can hold, and how its initial conditions could or should be
chosen.

It is an obvious idea to choose the initial conditions according to
$T_{-1}^{(n)} = 0$ and $T_{0}^{(n)} = S_{n}$ as in the case of
epsilon. However, these assumptions do not suffice here. Setting $k=-1$
in \cref{Diff_Eq_1stOrd_k_Tkn_1} yields:
\begin{equation}
  \label{Diff_Eq_1stOrd_k_Tkn_2}
  T_{1}^{(n)} \; = \; T_{-1}^{(n+1)} + \frac
   {\left[ T_{0}^{(n)} - T_{-2}^{(n+1)} \right] \,
    \left[ T_{-1}^{(n+1)} - T_{-1}^{(n)} \right]}
     {T_{0}^{(n+1)} - T_{0}^{(n)}} \, .
\end{equation}
The numerator on the right-hand side contains the undefined quantity
$T_{-2}^{(n+1)}$ and the undefined finite difference
$\Delta T_{-1}^{(n)} = T_{-1}^{(n+1)} - T_{-1}^{(n)}$. To remove these
ambiguities, we could for example assume
$\Delta T_{-1}^{(n)} = T_{-1}^{(n+1)} - T_{-1}^{(n)} = 0$ and
$T_{-2}^{(n+1)}[T_{-1}^{(n+1)} - T_{-1}^{(n)}] = -1$. But in our opinion
it would be simpler and more natural to use the recursion
(\ref{Diff_Eq_1stOrd_k_Tkn_1}) in combination with suitable initial
conditions only for $k \ge 0$.

If we set $k=0$ in \cref{Diff_Eq_1stOrd_k_Tkn_1} and use the epsilon
initial conditions $T_{-1}^{(n)} = 0$, $T_{0}^{(n)} = S_{n}$, and
$T_{1}^{(n)} = 1/[T_{0}^{(n+1)} - T_{0}^{(n)}]$, we obtain
\begin{equation}
  \label{Diff_Eq_1stOrd_k_Tkn_4}
  T_{2}^{(n)} \; = \;
   T_{0}^{(n+1)} + \frac{1}{T_{1}^{(n+1)} - T_{1}^{(n)}} \, .
\end{equation}
Comparison with \cref{eps_al_b} shows that $T_{2}^{(n)}$ is nothing but
$\varepsilon_{2}^{(n)}$ in disguise. Next, we set $k=1$ in
\cref{Diff_Eq_1stOrd_k_Tkn_1}:
\begin{equation}
  \label{Diff_Eq_1stOrd_k_Tkn_5}
  T_{3}^{(n)} \; = \; T_{1}^{(n+1)} + \frac
   {\left[ T_{2}^{(n)} - T_{0}^{(n+1)} \right] \,
    \left[ T_{1}^{(n+1)} - T_{1}^{(n)} \right]}
     {T_{2}^{(n+1)} - T_{2}^{(n)}} \, .
\end{equation}
At first sight, this looks like a new expression
$T_{3}^{(n)} \ne \varepsilon_{3}^{(n)}$. Unfortunately,
\cref{Diff_Eq_1stOrd_k_Tkn_4} implies
$T_{2}^{(n)} - T_{0}^{(n+1)} = 1/[ T_{1}^{(n+1)} - T_{1}^{(n)}]$. Thus,
we obtain $T_{3}^{(n)} = \varepsilon_{3}^{(n)}$, and with the help of
complete induction in $k$ we deduce $T_{k}^{(n)} = \varepsilon_{k}^{(n)}$
for $k \ge 3$, or actually for all $k \in \mathbb{N}_{0}$. These
considerations show that the recursion \eqref{Diff_Eq_1stOrd_k_Tkn_1}
produces nothing new if we use it in combination with the epsilon initial
conditions $T_{-1}^{(n)} = 0$, $T_{0}^{(n)} = S_{n}$, and
$T_{1}^{(n)} = 1/[T_{0}^{(n+1)} - T_{0}^{(n)}]$.

We will now discuss the choice of \emph{alternative} initial conditions
for the recursion \eqref{Rec_Tkn_Fn0_Fn1}, which is based on the
assumption that $F_{k}^{(n)}$ satisfies the second order finite
difference equation \eqref{Diff_Eq_2ndOrd_n-dep_k}. This is a crucial
step, and it is also the point where our approach becomes experimental.

We choose parts of the initial values of the recursion
\eqref{Diff_Eq_2ndOrd_k_Tkn} as in Wynn's epsilon algorithm:
\begin{equation}
  \label{EpsIniVa}
  T_{-1}^{(n)} \; = \; 0 \, ,
   \quad T_{0}^{(n)} \; = \; S_{n} \, ,
    \quad T_{1}^{(n)} \; = \; \frac {1}{\Delta S_{n}} \, ,
     \qquad n \in \mathbb{N}_{0} \, .
\end{equation}
With these initial conditions, epsilon and rho are recovered by choosing
$T_{2}^{(n)} = \varepsilon_{2}^{(n)}$ and $T_{2}^{(n)} = \rho_{2}^{(n)}$,
respectively.

If we want to employ the recursion \eqref{Diff_Eq_2ndOrd_k_Tkn}, we first
have to make a choice about the initial value
$F_{1}^{(n)}$. \Cref{Tkn_2_Fkn} implies
$F_{1}^{(n)} = [T_{2}^{(n)} - T_{0}^{(n+1)}]
[T_{1}^{(n+1)}-T_{1}^{(n)}]$.
If $T_{-1}^{(n)}$, $T_{0}^{(n)}$, and $T_{1}^{(n)}$ are chosen
 according
to \cref{EpsIniVa}, we only need one further choice for $T_{2}^{(n)}$ to
fix $F_{1}^{(n)}$.

It is a relatively natural idea to choose $T_{2}^{(n)}$ in such a way
that it corresponds to the special case of a sequence transformation with
low transformation order. If we want the resulting new transformation to
be more versatile than either epsilon or rho, it makes sense to choose a
transformation which is known to accelerate both linear and logarithmic
convergence effectively. Obvious candidates would be the spacial case
$\vartheta_{2}^{(n)}$ of Brezinski's theta algorithm discussed in
\cref{Sub:BrezThetaAlg}, or suitable variants of Levin's transformation
discussed in \cref{Sub:LevinU,Sub:LevinV}.

There is another decision, which we have to make. Should the expression
for $T_{2}^{(n)}$ depend on $n$ only \emph{implicitly} via the input data
$\{ S_{n} \}$, or can it depend on $n$ also \emph{explicitly} (examples
of low order transformations of that kind can be found in \citep[Section
3]{Weniger/1991}). We think that it would be more in agreement with the
spirit of the other initial conditions \eqref{EpsIniVa} if we require
that $T_{2}^{(n)}$ depends on $n$ only \emph{implicitly} via the input
data.

Under these assumptions, it looks like an obvious idea to identify
$T_{2}^{(n)}$ either with Brezinski's theta transformation
$\vartheta_{2}^{(n)}$ given by \cref{theta_2_1} or with the explicit
expressions for the Levin $u$ and $v$ transformations
$u_{2}^{(n)}(\beta, S_{n})$ and $v_{1}^{(n)} (\beta, S_{n})$ given by
\cref{theta_2<->u_2,theta_2<->v_1}. 

We checked the performance of the resulting algorithms numerically. For
$T_{2}^{(n)} = \vartheta_{2}^{(n)} = u_{2}^{(n+1)} (\beta,S_{n+1}) =
v_{1}^{(n+1)}(\beta, S_{n+1})$ (compare
\cref{theta_2<->u_2,theta_2<->v_1}), we obtained consistently good
results, but our results for
$T_{2}^{(n)} = u_{2}^{(n)} (\beta,S_{n}) = v_{1}^{(n)}(\beta, S_{n})$
(compare \cref{uLevTr_5,vLevTr_2}) were less good. Thus, we will
exclusively use the following non-epsilon initial condition:
\begin{align}
  \label{IniConds}
  T_{2}^{(n)} & \; = \; \vartheta_{2}^{(n)} \; = \;
   u_{2}^{(n+1)} (\beta,S_{n+1}) \; = \;
    v_{1}^{(n+1)} (\beta, S_{n+1})
  \notag
  \\              
  & \; = \; S_{n+1} - \frac
   {[\Delta^{2} S_{n+1}] [\Delta S_{n}] [\Delta S_{n+1}]}
    {[\Delta S_{n+2}] [\Delta^{2} S_{n}] -
     [\Delta S_{n}] [\Delta^{2} S_{n+1}]} \, ,
      \qquad n \in \mathbb{N}_{0} \, .
\end{align}
Then, \cref{Tkn_2_Fkn} implies
\begin{equation}
  \label{Choice_F0F1}
  F_{0}^{(n)} \; = \; 1 \, ,
   \quad F_{1}^{(n)} \; = \;  \frac
    {\bigl[ \Delta^{2} T_{0}^{(n)} \bigr]
     \bigl[ \Delta^{2} T_{0}^{(n+1)} \bigr]}
    {\bigl[ \Delta T_{0}^{(n+2)} \bigr]
     \bigl[\Delta^{2} T_{0}^{(n)} \bigr]
      - \bigl[\Delta  T_{0}^{(n)}\bigr]
       \bigl[\Delta^{2} T_{0}^{(n+1)} \bigr]}
\end{equation}
and we obtain the following new recursive scheme:
\begin{subequations}
  \label{seps}
  \begin{align}
    \label{seps_a}
    \tilde{T}_{0}^{(n)} & \; = \; S_{n} \, ,
    \\
    \label{seps_b}
    \tilde{T}_{1}^{(n)} & \; = \; \frac
     {1} {\tilde{T}_{0}^{(n+1)} - \tilde{T}_{0}^{(n)}} \, ,
    \\
    \label{seps_c}
    \tilde{T}_{2}^{(n)} & \; = \; \tilde{T}_{0}^{(n+1)} - \frac
     {\bigl[ \Delta \tilde{T}_{0}^{(n)} \bigr]
      \bigl[ \Delta \tilde{T}_{0}^{(n+1)} \bigr]
       \bigl[ \Delta^{2} \tilde{T}_{0}^{(n+1)} \bigr]}
        {\bigl[ \Delta \tilde{T}_{0}^{(n+2)} \bigr]
         \bigl[ \Delta^{2} \tilde{T}_{0}^{(n)} \bigr] -
          \bigl[ \Delta \tilde{T}_{0}^{(n)} \bigr]
           \bigl[ \Delta^{2} \tilde{T}_{0}^{(n+1)} \bigr]} \, ,
    \\
    \label{seps_d}
    F_{1}^{(n)} & \; = \; \frac
     {\bigl[ \Delta^{2} \tilde{T}_{0}^{(n)} \bigr]
      \bigl[ \Delta^{2} \tilde{T}_{0}^{(n+1)} \bigr]}
     {\bigl[ \Delta \tilde{T}_{0}^{(n+2)} \bigr]
      \bigl[ \Delta^{2} \tilde{T}_{0}^{(n)} \bigr] -
       \bigl[ \Delta \tilde{T}_{0}^{(n)} \bigr]
        \bigl[ \Delta^{2} \tilde{T}_{0}^{(n+1)} \bigr]} \, ,
   \\
    \label{seps_e}
    \tilde{T}_{k+1}^{(n)} & \; = \; \tilde{T}_{k-1}^{(n+1)} \, + \, \frac
     {1-k+k F_{1}^{(n)}} {\tilde{T}_{k}^{(n+1)} - \tilde{T}_{k}^{(n)}}
      \, , \qquad k = 2, 3, \ldots \, ,
         \qquad n \in \mathbb{N}_{0} \, .
  \end{align}
\end{subequations}
In \cref{seps}, we wrote $\tilde{T}$ instead of $T$ in order to
distinguish our new algorithm from others mentioned before. The
performance of our new algorithm \eqref{seps} will be studied in detail
in the following \cref{Sec:NumEx}.

%
\typeout{==> Section: Numerical examples}
\section{Numerical examples}
\label{Sec:NumEx}
%

In this section, we apply our new algorithm \eqref{seps} to some linearly
and logarithmically convergent sequences and also to some divergent
series. We compare the numerical results of our new algorithm with those
obtained by applying Wynn's epsilon and rho algorithms \eqref{eps_al} and
\eqref{RhoAlStand}, Osada's generalized rho algorithm \eqref{OsRhoAl},
Brezinski's theta algorithm \eqref{ThetaAl} and the iteration
\eqref{ThetIt} of $\vartheta_{2}^{(n)}$. All the numerical computation are
obtained with the help of MATLAB 7.0 in double precision (15 decimal digits).

\typeout{==> Subsection: Asymptotic model sequences}
\subsection{Asymptotic model sequences}
\label{Sub:AsyModSeq}
%

In this Section, we want to analyze the numerical performance of our new
sequence transformation \eqref{seps}. As input data, we use some model
sequences which provide reasonably good approximations to the elements of
numerous practically relevant sequences in the limit of large
indices. Essentially the same model sequences had been used by Sidi in
\citep[Eqs.\ (1.1) - (1.3)]{Sidi/2002} and in \citep[Definition 15.3.2 on
p.\ 285]{Sidi/2003}, and by \citet[Section IV]{Weniger/2004}:
\begin{enumerate}
\item Logarithmically convergent model sequence
  \begin{equation}
    \label{eq:log}
    S_{n} \; = \; S + \sum_{j=0}^{\infty} \, c_{j} \, n^{\eta-j} \, ,
     \qquad c_{0} \neq 0 \, , \quad n \in \mathbb{N}_{0} \, .
  \end{equation}
\item Linearly convergent model sequence
  \begin{equation}
    \label{eq:linear}
    S_{n} \; = \; S + \xi^{n} \, \sum_{j=0}^{\infty} \,
     c_{j} \, n^{\eta-j} \, , \qquad c_{0} \neq 0 \, ,
      \quad n \in \mathbb{N}_{0} \, .
  \end{equation}
\item Hyperlinearly convergent model sequence
  \begin{equation}
    \label{eq:fac}
    S_{n} \; = \; S + \frac{\xi^{n}}{(n!)^{r}} \,
     \sum_{j=0}^{\infty} \, c_{j} \, n^{\eta-j} \, ,
      \qquad c_{0} \neq 0 \, , \quad r \in \mathbb{R} \, ,
       \quad n \in \mathbb{N}_{0} \, .
  \end{equation}
\end{enumerate}
The series expansions \crefrange{eq:log}{eq:fac} are to be interpreted as
purely \emph{formal} expansions, i.e., we do not tacitly assume that they
converge. It is actually more realistic to assume that they are only
asymptotic to $S_{n}$ as $n \to \infty$. There are, however, restrictions
on the parameters $\eta$ and $\xi$ if we require that the sequences
$\left\{ S_{n} \right\}$ defined by the series expansions
\crefrange{eq:log}{eq:fac} are able to represent something finite.

The sequence $\{ S_{n} \}$ defined by \cref{eq:log} can only converge to
its limit $S$ if $\eta < 0$ holds. For $\eta > 0$, \citet[Text following
Eq.\ (1.1)]{Sidi/2002} calls $S$ the \emph{antilimit} of this
sequence. While this terminology is very useful in the case of divergent,
but summable series, it is highly questionable here. The divergence of
this sequence for $\eta > 0$ is genuine, and the normally used summation
techniques are not able to produce anything finite. Thus, the terminology
\emph{antilimit} does not make sense here.

The sequence $\left\{ S_{n} \right\}$ defined by \cref{eq:linear}
converges for all $\eta \in \mathbb{R}$ linearly to its limit $S$
provided that $\vert \xi \vert < 1$ holds. For $\vert \xi \vert > 1$,
this sequence diverges, but it should normally be summable to $S$ for all
$\xi$ belonging to the cut complex plane
$\mathbb{C} \setminus [1, \infty)$.

The most convenient situation occurs if the sequence elements $S_{n}$
satisfy \cref{eq:fac}. Then, this sequence converges for all
$\xi, \eta \in \mathbb{R}$ hyperlinearly to its limit $S$. In
\cref{eq:fac}, the restriction $r \in \mathbb{R}_+$ is
essential, because for $r<0$, we
would obtain a sequence which diverges (hyper)factorially for all
$\xi \ne 0$.  Divergent sequences of that kind can according to
experience often be summed to their generalized limits $S$ for all
$\xi \in \mathbb{C} \setminus [0, \infty)$, but a rigorous theoretical
analysis of such a summation process is normally extremely difficult
(compare the rigorous theoretical analysis of the summation of a
comparatively simple special case in \citep{Borghi/Weniger/2015} and
references therein).

%
\typeout{==> Subsection: Linear convergence}
\subsection{Linearly convergent sequences}
\label{SubSec:LinConv}
%

In this Subsection, we apply our new algorithm to some linearly
convergent alternating and monotone series, which are special cases of
the model sequence \eqref{eq:log}. The numerical results obtained by our
new algorithm \eqref{seps}, Wynn's epsilon algorithm \eqref{eps_al},
Brezinski's theta algorithm \eqref{ThetaAl} and the iteration
\eqref{ThetIt} of $\vartheta_{2}^{(n)}$ are presented in
\cref{linear1,linear2}.

In the captions of the following Tables, $\Ent{x}$ denotes the
\emph{integral part} of $x$, which is the \emph{largest} integer $n$
satisfying $n \le x$.

\textbf{Example 4.1:} We consider the alternating partial sums
\begin{equation}
  \label{lisq1}
  S_{n} \; = \; \sum_{k=0}^{n} \frac{(-1)^{k}}{k+1} \, ,
   \qquad n \in \mathbb{N}_{0} \, ,
\end{equation}
which converge to $S=\ln (2)$. The \emph{Boole summation formula}
\citep*[Eq.\ (24.17.1)]{Olver/Lozier/Boisvert/Clark/2010} yields
\begin{equation}
  S_{n} - \ln (2) \; \sim \; \frac{(-1)^{n-1}}{2}
   \left\{ \frac {1} {n+1} + \mathrm{O} \bigl(n^{-2} \bigr) \right\}
    \, , \qquad n \to \infty \, .
\end{equation}
This truncated asymptotic expansion is a special case of the model
sequence \eqref{eq:linear} with $\xi=-1$ and $\eta=-1$.

\begin{table}
\caption{\label{linear1}Numerical results of example 4.1 }
\begin{tabular}{|c|c|c|c|c|c|}\hline
$n$ & $ \vert S_n-S \vert$ &
$\bigl\vert \tilde{T}_{2\Ent{(n-1)/2}}^{(n-1-2\Ent{(n-1)/2})}-S \bigr\vert$ &
$\bigl\vert \epsilon_{2\Ent{n/2}}^{(n-2\Ent{n/2})}-S \bigr\vert$
& $\bigl\vert \vartheta_{2\Ent{n/3}}^{(n-3\Ent{n/3})}-S \bigr\vert$ &
$\bigl\vert \mathcal{J}_{\Ent{n/3}}^{(n-3\Ent{n/3})}-S \bigr\vert$ \\
&\cref{lisq1} & \cref{seps} & \cref{eps_al} & \cref{ThetaAl} &
\cref{ThetIt} \\ \hline
7 &  0.059 & 6.768$\times10^{-7}$  & 1.437$\times10^{-6}$  & 5.702$\times10^{-7}$  & 1.527$\times10^{-8}$ \\
8 &  0.052 & 3.485$\times10^{-8}$  & 1.518$\times10^{-7}$  & 1.961$\times10^{-7}$  & 7.793$\times10^{-9}$ \\
9 & 0.048  & 5.907$\times10^{-9}$  & 3.807$\times10^{-8}$  & 8.881$\times10^{-10}$ & 6.034$\times10^{-10}$ \\
10 & 0.043 & 1.920$\times10^{-9}$ & 4.402$\times10^{-9}$   & 2.764$\times10^{-10}$ & 2.115$\times10^{-10}$ \\
11 & 0.040 & 8.728$\times10^{-11}$ & 1.042$\times10^{-9}$  & 9.280$\times10^{-11}$ & 7.530$\times10^{-11}$ \\
12 & 0.037 & 4.244$\times10^{-11}$ & 1.282$\times10^{-10}$ & 1.171$\times10^{-13}$ & 4.954$\times10^{-13}$ \\
13 & 0.034 & 1.396$\times10^{-11}$ & 2.909$\times10^{-11}$ & 1.665$\times10^{-14}$ & 3.331$\times10^{-14}$ \\
14 & 0.032 & 3.936$\times10^{-13}$ & 3.745$\times10^{-12}$ & 1.266$\times10^{-14}$ & 2.331$\times10^{-15}$ \\
\hline
\end{tabular}
\end{table}

\textbf{Example 4.2:} The linearly convergent monotone sequence
\begin{equation}
  \label{lisq2}
  S_{n} \; = \; \frac{4}{5} \, \sum_{k=0}^{n} \,
   \frac{(4/5)^{k}}{k+1} \, ,
    \qquad n \in \mathbb{N}_{0} \, ,
\end{equation}
converges to $S = \ln (5)$ as $n \to \infty$. As shown by \citet[p.\
84]{Sidi/2003}, $S_{n}$ possesses the asymptotic expansion
\begin{equation}
  S_{n} - \ln(5) \; \sim \; \frac{(4/5)^{n}}{n} \,
   \left\{ -4 + \mathrm{O} (n^{-2}) \right\} \, ,
    \qquad n \to \infty \, ,
\end{equation}
which is a special case of the model sequence \eqref{eq:linear} with
$\xi=4/5$ and $\eta=-1$.

\begin{table}
\caption{\label{linear2}Numerical results of example 4.2 }
\begin{tabular}{|c|c|c|c|c|c|}\hline
$n$ & $\vert S_n-S \vert$ & $\bigl\vert
\tilde{T}_{2\Ent{(n-1)/2}}^{(n-1-2\Ent{(n-1)/2})}-S  \bigr\vert$
& $\bigl\vert \epsilon_{2\Ent{n/2}}^{(n-2\Ent{n/2})}-S \bigr\vert$ &
$\bigl\vert \vartheta_{2\Ent{n/3}}^{(n-3\Ent{n/3})}-S \bigr\vert$
& $\bigl\vert \mathcal{J}_{\Ent{n/3}}^{(n-3\Ent{n/3})}-S \bigr\vert$ \\
& \cref{lisq2} & \cref{seps} & \cref{eps_al} & \cref{ThetaAl} &
\cref{ThetIt} \\ \hline
13 & 9.657$\times10^{-3}$ & 7.687$\times10^{-4}$  & 3.637$\times10^{-6}$  & 3.533$\times10^{-6}$ & 9.972$\times10^{-10}$ \\
14 & 7.312$\times10^{-3}$ & 1.226$\times10^{-4}$  & 1.272$\times10^{-6}$  & 3.560$\times10^{-6}$ & 2.092$\times10^{-10}$ \\
15 & 5.552$\times10^{-3}$ & 3.558$\times10^{-4}$  & 5.260$\times10^{-7}$  & 3.561$\times10^{-6}$ & 2.026$\times10^{-9}$ \\
16 & 4.228$\times10^{-3}$ & 4.880$\times10^{-5}$  & 1.860$\times10^{-7}$  & 3.585$\times10^{-6}$ & 2.479$\times10^{-10}$ \\
17 & 3.227$\times10^{-3}$ & 1.608$\times10^{-5}$  & 7.597$\times10^{-8}$  & 3.573$\times10^{-6}$ & 6.144$\times10^{-11}$ \\
18 & 2.469$\times10^{-3}$ & 1.990$\times10^{-5}$  & 2.769$\times10^{-8}$  & 3.582$\times10^{-6}$ & 3.903$\times10^{-11}$ \\
19 & 1.892$\times10^{-3}$ & 8.434$\times10^{-6}$  & 1.025$\times10^{-8}$  & 1.837$\times10^{-6}$ & 5.445$\times10^{-11}$ \\
20 & 1.453$\times10^{-3}$ & 4.816$\times10^{-6}$  & 5.256$\times10^{-9}$  & 5.225$\times10^{-9}$ & 4.459$\times10^{-11}$ \\

\hline
\end{tabular}
\end{table}

\typeout{==> Subsection: Logarithmic convergence}
\subsection{Logarithmically convergent sequences}
\label{SubSec:LogConv}

First, we consider a logarithmically convergent sequence of the type of
the model sequence \eqref{ModSeqAlpha} with integral decay parameter
$\theta$. In Table \ref{logsq1}, we present the numerical results
obtained by applying our new algorithm \eqref{seps}, the standard form
\eqref{RhoAlStand} of Wynn's rho algorithm, Brezinski's theta algorithm
\eqref{ThetaAl}, and the iteration \eqref{ThetIt} of
$\vartheta_{2}^{(n)}$. Other logarithmically convergent sequences with
non-integral decay parameters $\theta$ will be discussed later.

\textbf{Example 4.3:} We consider the logarithmically convergent sequence
\begin{equation}\label{logex1}
 S_n \; = \; 1+ \sum_{k=1}^n
   \left\{ \frac{1}{k+1}+\log \left(\frac{k}{k+1} \right) \right\} \, ,
   \qquad n \in \mathbb{N} \, ,
\end{equation}
which converges to
$S=0.57721~56649~01532 \cdots$ (Euler's constant).

The Euler-Maclaurin formula (see for example \citep*[Eq.\
(2.10.1)]{Olver/Lozier/Boisvert/Clark/2010}) yields:
\begin{equation}
  S_{n-1} \; = \;
   \sum_{k=1}^{n} \frac{1}{k} - \log (n) \; \sim \;
    S + n^{-1} \left( 1/2 + \sum_{j=1}^{\infty}
     \frac{B_{2j}}{2j} n^{-2j} \right) \, , \qquad n \to \infty \, ,
\end{equation}
where $B_{n}$ is a Bernoulli number \citep*[Eq.\ (24.2.1)]{Olver/Lozier/Boisvert/Clark/2010}.

\begin{table}
\caption{\label{logsq1}Numerical results of example 4.3 }
\begin{tabular}{|c|c|c|c|c|c|}\hline
$n$ & $\bigl\vert S_n-S \bigr\vert$ & $\bigl\vert
\tilde{T}_{2\Ent{(n-1)/2}}^{(n-1-2\Ent{(n-1)/2})}-S\bigr\vert$ &
$\bigl\vert \rho_{2\Ent{n/2}}^{(n-2\Ent{n/2})}-S \bigr\vert$ &
$\bigl\vert \vartheta_{2\Ent{n/3}}^{(n-3\Ent{n/3})}-S \bigr\vert$ &
$\bigl\vert \mathcal{J}_{\Ent{n/3}}^{(n-3\Ent{n/3})}-S \bigr\vert$ \\
&\cref{{logex1}} &\cref{seps} &\cref{RhoAlStand} &\cref{ThetaAl}
&\cref{ThetIt} \\ \hline
3 &  0.120 & 4.051$\times10^{-4}$  & 2.942$\times10^{-4}$  & 4.051$\times10^{-4}$  & 4.051$\times10^{-4}$ \\
4 &  0.097 & 3.190$\times10^{-4}$  & 6.774$\times10^{-5}$  & 3.190$\times10^{-4}$  & 3.190$\times10^{-4}$ \\
5 &  0.081 & 6.518$\times10^{-4}$  & 7.156$\times10^{-6}$  & 1.875$\times10^{-4}$  & 1.875$\times10^{-4}$ \\
6 &  0.070 & 2.724$\times10^{-5}$  & 8.286$\times10^{-8}$  & 2.391$\times10^{-4}$  & 2.393$\times10^{-4}$ \\
7 &  0.061 & 1.521$\times10^{-8}$  & 5.091$\times10^{-9}$  & 2.221$\times10^{-5}$  & 2.238$\times10^{-5}$ \\
8 &  0.055 & 1.188$\times10^{-8}$  & 1.870$\times10^{-9}$  & 4.735$\times10^{-6}$  & 4.822$\times10^{-6}$ \\
9 &  0.049 & 5.295$\times10^{-9}$  & 6.288$\times10^{-10}$ & 2.316$\times10^{-7}$  & 1.974$\times10^{-7}$ \\
10 & 0.045 & 6.931$\times10^{-10}$ & 1.289$\times10^{-11}$ & 4.116$\times10^{-8}$  & 7.269$\times10^{-9}$ \\
\hline
\end{tabular}
\end{table}

Next, we consider logarithmically convergent sequences of the form of
\cref{eq:log} with non-integral decay $\theta$. In Table \ref{logsq2}, we
present the numerical results for the test problem \eqref{LemConst}
obtained by our new algorithm \eqref{seps}, Osada's generalized rho
algorithm \eqref{OsRhoAl} with $\theta=1/2$, Brezinski's theta algorithm
\eqref{ThetaAl} and the iteration \eqref{ThetIt} of
$\vartheta_{2}^{(n)}$.

\textbf{Example 4.4:} A very useful test problem is the series expansion
for the so-called lemniscate constant $A$ \citep[Theorem 5]{Todd/1975}:
\begin{equation}
  \label{LemConst}
  A \; = \; \frac{[\Gamma (1/4)]^2}{4 \, (2 \pi)^{1/2}} \; = \;
   \sum_{m=0}^{\infty} \,
    \frac{(2m-1)!!}{(2m)!!} \, \frac{1}{4m+1} \, .
\end{equation}
Standard asymptotic techniques show that the terms
$(2m-1)!!/[(2m)!! (4m+1)]$ of this series decay like $m^{-3/2}$ as
$m \to \infty$. This implies that the remainders
\begin{equation}
  \label{ParSum_lem}
  R_{n} \; = \; S_{n} - S \; = \;
   \sum_{m=n+1}^{\infty} \, \frac{(2m-1)!!}{(2m)!!} \, \frac{1}{4m+1}
\end{equation}
decay like $n^{-1/2}$ as $n \to \infty$. Consequently, the series
expansion \eqref{LemConst} constitutes an extremely challenging
convergence acceleration problem. \citet[p.\ 14]{Todd/1975} stated in
\citeyear{Todd/1975} that this series expansion is of no practical
use for the computation of $A$, which is certainly true if no convergence
acceleration techniques are available. In \citeyear{Smith/Ford/1979}, the
series expansion \eqref{LemConst} was used by \citet[Entry 13 of Table
6.1]{Smith/Ford/1979} to test the ability of sequence transformations of
accelerating logarithmic convergence. They observed that the standard
form \eqref{RhoAlStand} of Wynn's rho algorithm fails to accelerate the
convergence of this series \citep[p.\ 235]{Smith/Ford/1979} which is in
agreements with Osada's later convergence analysis \citep[Theorem
3.2]{Osada/1990a}. The series expansion \eqref{LemConst} was also
extensively used by \citet[Tables 14-2 - 14-4]{Weniger/1989} as a test
case.

\begin{table}
\caption{\label{logsq2}Numerical results of example 4.4 }
\begin{tabular}{|c|c|c|c|c|c|}\hline
$n$ & $\bigl\vert S_n-S \bigr\vert$ & $\bigl\vert
\tilde{T}_{2\Ent{(n-1)/2}}^{(n-1-2\Ent{(n-1)/2})}-S \bigr\vert$ &
$\bigl\vert {\bar \rho}_{2\Ent{n/2}}^{(n-2\Ent{n/2})}-S \bigr\vert$ &
$\bigl\vert \vartheta_{2\Ent{n/3}}^{(n-3\Ent{n/3})}-S \bigr\vert$ &
$\bigl\vert \mathcal{J}_{\Ent{n/3}}^{(n-3\Ent{n/3})}-S \big\vert$ \\
&\cref{ParSum_lem} & \cref{seps} & \cref{OsRhoAl}, $\theta=1/2$
&\cref{ThetaAl} & \cref{ThetIt} \\ \hline
14 & 0.073 & 1.029$\times10^{-8}$   & 8.637$\times10^{-12}$  & 8.733$\times10^{-10}$  & 2.680$\times10^{-9}$ \\
15 & 0.071 & 1.201$\times10^{-10}$  & 1.038$\times10^{-12}$  & 5.325$\times10^{-10}$  & 2.565$\times10^{-9}$ \\
16 & 0.069 & 5.704$\times10^{-11}$  & 2.517$\times10^{-12}$  & 5.542$\times10^{-10}$  & 1.921$\times10^{-9}$ \\
17 & 0.067 & 1.081$\times10^{-10}$  & 9.550$\times10^{-13}$  & 1.299$\times10^{-9}$   & 9.512$\times10^{-9}$ \\
18 & 0.065 & 1.606$\times10^{-11}$  & 1.384$\times10^{-11}$  & 5.174$\times10^{-10}$  & 2.457$\times10^{-9}$ \\
19 & 0.063 & 1.333$\times10^{-10}$  & 6.473$\times10^{-13}$  & 1.086$\times10^{-9}$   & 7.667$\times10^{-9}$ \\
20 & 0.062 & 1.583$\times10^{-11}$  & 3.970$\times10^{-13}$  & 1.614$\times10^{-6}$   & 7.311$\times10^{-8}$ \\
21 & 0.060 & 1.591$\times10^{-12}$  & 6.120$\times10^{-13}$  & 2.152$\times10^{-7}$   & 1.540$\times10^{-8}$ \\
\hline
\end{tabular}
\end{table}

\textbf{Example 4.5:} The series expansion \eqref{1/z_RBF} of $1/z$ in
terms of reduced Bessel functions, whose basic properties are reviewed in
\cref{App:RBF}, is a challenging test problem for the ability of a
sequence transformation to accelerate the convergence of a
logarithmically convergent sequence with a non-integral decay parameter.

We consider the sequence of partial sum of the series \eqref{1/z_RBF}
with $z=1$,
\begin{equation}
  \label{besfsq}
  S_{n} \; = \; \sum_{m=0}^{n} \hat{k}_{m-1/2} (1) /[2^m \, m!],
\end{equation}
which corresponds to a special case of the model sequence \eqref{eq:log}
with $\eta =- 1/2$. In \cref{logsq3}, we present the numerical results
for the test problem \eqref{besfsq} obtained by our new algorithm
\eqref{seps}, Osada's generalized rho algorithm \eqref{OsRhoAl} with
$\theta=1/2$, Brezinski's theta algorithm \eqref{ThetaAl} and the
iteration \eqref{ThetIt} of $\vartheta_{2}^{(n)}$.

\begin{table}
\caption{\label{logsq3}Numerical results of example 4.5 }
\begin{tabular}{|c|c|c|c|c|c|}\hline
$n$ & $\bigl\vert S_n-S \bigr\vert$ &
$\bigl\vert \tilde{T}_{2\Ent{(n-1)/2}}^{(n-1-2\Ent{(n-1)/2})}-S \bigr\vert$ &
$\bigl\vert {\bar \rho}_{2\Ent{n/2}}^{(n-2\Ent{n/2})}-S \bigr\vert$ &
$\bigl\vert \vartheta_{2\Ent{n/3}}^{(n-3\Ent{n/3})}-S \bigr\vert$ &
$\bigl\vert
\mathcal{J}_{\Ent{n/3}}^{(n-3\Ent{n/3})}-S \bigr\vert$ \\
&\cref{besfsq} & \cref{seps} & \cref{OsRhoAl}, $\theta=1/2$
&\cref{ThetaAl} &\cref{ThetIt} \\ \hline
14 & 0.149 & 4.688$\times10^{-9}$   & 3.580$\times10^{-10}$  & 4.073$\times10^{-5}$ & 1.149$\times10^{-7}$ \\
15 & 0.144 & 3.321$\times10^{-9}$   & 1.266$\times10^{-10}$  & 4.362$\times10^{-5}$ & 1.068$\times10^{-8}$ \\
16 & 0.139 & 1.443$\times10^{-8}$   & 2.067$\times10^{-10}$  & 2.654$\times10^{-5}$ & 8.333$\times10^{-10}$ \\
17 & 0.135 & 1.086$\times10^{-9}$   & 1.338$\times10^{-10}$  & 5.316$\times10^{-8}$ & 4.021$\times10^{-9}$ \\
18 & 0.131 & 3.993$\times10^{-9}$   & 3.329$\times10^{-11}$  & 8.259$\times10^{-8}$ & 2.607$\times10^{-9}$ \\
19 & 0.128 & 4.042$\times10^{-10}$  & 8.087$\times10^{-11}$  & 5.279$\times10^{-8}$ & 7.819$\times10^{-10}$ \\
20 & 0.125 & 2.013$\times10^{-9}$   & 6.629$\times10^{-11}$  & 5.270$\times10^{-8}$ & 1.851$\times10^{-8}$ \\
21 & 0.122 & 1.075$\times10^{-11}$  & 8.181$\times10^{-11}$  & 5.275$\times10^{-8}$ & 1.858$\times10^{-8}$ \\
\hline
\end{tabular}
\end{table}

Our numerical Examples 4.3 - 4.5 show that our new algorithm is
apparently able to accelerate the convergence of logarithmically
convergent sequences of the type of \eqref{ModSeqAlpha} \emph{without}
having to know the value of the decay parameter $\theta$ explicitly, as
it is necessary in the case of Osada's variant \eqref{OsRhoAl} of the rho
algorithm.

\typeout{==> Subsection: Divergent series}
\subsection{Divergent series}
\label{SubSec:DivSer}

In this Section, we apply our new algorithm for the summation of the
factorially divergent Euler series. Since the early days of calculus,
divergent series have been a highly controversial topic in mathematics
\citep{Burkhardt/1911,Ferraro/2008,Tucciarone/1973}, and to some extend
they still are (see for example \citep[Appendix D]{Weniger/2008} and
references therein).

The exponential integral \citep*[Eq.\
(6.2.1)]{Olver/Lozier/Boisvert/Clark/2010}
\begin{equation}
  \label{Def_ExpInt}
  E_1 (z) \; = \;
  \int_{z}^{\infty} \, \frac{\exp (-t)}{t} \, \mathrm{d} t \, ,
\end{equation}
possesses the asymptotic expansion \citep*[Eq.\
(6.12.1)]{Olver/Lozier/Boisvert/Clark/2010}
\begin{equation}
  \label{ExpInt_AsySer}
  z \, \mathrm{e}^z \, E_1 (z) \; \sim \;
  \sum_{m=0}^{\infty} \, \frac{(-1)^m m!}{z^m}
  \; = \; {}_2 F_0 (1, 1; -1/z) \, , \qquad z \to \infty \, .
\end{equation}
The exponential integral can be expressed as a Stieltjes integral
\citep*[Eq.\ (6.7.1 )]{Olver/Lozier/Boisvert/Clark/2010},
\begin{equation}
  \label{ExpInt_StieInt}
  \mathrm{e}^z \, E_1 (z) \; = \;
  \int_{0}^{\infty} \, \frac{\exp (-t) \mathrm{d}t}{z+t} \; = \;
  \frac{1}{z} \,
  \int_{0}^{\infty} \, \frac{\exp (-t) \mathrm{d}t}{1+t/z} \, ,
\end{equation}
which implies that the divergent inverse power series
(\ref{ExpInt_AsySer}) is a Stieltjes series.

Asymptotic series as $z \to \infty$ involving a factorially divergent
generalized hypergeometric series ${}_2 F_0$ occur also among other
special functions. Examples are the asymptotic expansion of the modified
Bessel function $K_{\nu} (z)$ of the second kind \citep*[Eq.\
(10.40.2)]{Olver/Lozier/Boisvert/Clark/2010}, or the asymptotic expansion
of the Whittaker function $W_{\kappa,\mu} (z)$ of the second kind
\citep*[Eq.\ (13.19.3)]{Olver/Lozier/Boisvert/Clark/2010}. Similar
divergent series occur also in statistics, as discussed in a book by
\citet*{Bowman/Shenton/1989}.

If we replace in \cref{ExpInt_AsySer} $z$ by $1/z$, we obtain the
so-called Euler series, which had already been studied by Euler (see for
example the books by \citet[pp.\ 323 - 324]{Bromwich/1991} and
\citet[pp.\ 26 - 29]{Hardy/1949} or the articles by \citet{Barbeau/1979}
and \citet{Barbeau/Leah/1976}). For $\vert \arg(z) \vert < \pi$, the
Euler series is the asymptotic series of the Euler integral:
\begin{equation}
  \label{Def:EulerSer}
  \mathcal{E} (z) \; = \; \int_{0}^{\infty} \,
   \frac{\mathrm{e}^{-t}}{1 + zt} \, \mathrm{d} t \; \sim \;
    \sum _{\nu = 0}^{\infty} \; (-1)^{\nu} \, \nu! \, z^{\nu}
     \; = \; {}_2 F_0 (1,1;- z) \, , \qquad z \to 0 \, .
\end{equation}
The Euler series is the most simple prototype of the factorially
divergent perturbation series that occur abundantly in quantum
physics. Already \citet{Dyson/1952} had argued that perturbation
expansions in quantum electrodynamics must diverge factorially. Around
1970, \citet{Bender/Wu/1969,Bender/Wu/1971,Bender/Wu/1973} showed in
their seminal work on anharmonic oscillators that factorially divergent
perturbation expansions occur in nonrelativistic quantum
mechanics. Later, it was found that factorially divergent perturbation
expansions are actually the rule in quantum physics rather than the
exception (see for example the articles by \citet[Table 1]{Fischer/1997}
and \citet{Suslov/2005}, or the articles reprinted in the book by
\citet{LeGuillou/Zinn-Justin/1990}). Applications of factorially
divergent series in optics are discussed in a recent review by
\citet{Borghi/2016}.

The Euler integral $\mathcal{E} (z)$ is a Stieltjes function. Thus, the
corresponding Euler series ${}_2 F_0 (1,1;- z)$ is its associated
Stieltjes series (see for example \cite[Chapter
5]{Baker/Graves-Morris/1996}). This Stieltjes property guarantees the
convergence of certain subsequences of the Pad\'{e} table of
${}_2 F_0 (1,1;- z)$ to the Euler integral (see for example \cite[Chapter
5.2]{Baker/Graves-Morris/1996} or the discussion in
\citep*{Borghi/Weniger/2015}).

It has been shown in many calculations that sequence transformations are
principal tools that can accomplish an efficient summation of the
factorially divergent expansions of the type of the Euler series (see for
example \citep*{Borghi/Weniger/2015} and references therein). However, in
the case of most sequence transformations, no rigorous theoretical
convergence proofs are known (this applies also to Pad\'{e} approximants
if the series to be transformed is not Stieltjes). An exception is the
delta transformation \citep[Eq.\ (8.4-4)]{Weniger/1989} of the Euler
series, whose convergence to the Euler integral for all
$z \in \mathbb{C} \setminus [0, \infty)$ was proved rigorously in
\citep*{Borghi/Weniger/2015}.

In the following \cref{divsq1,divsq2}, we use as input data the partial sums
\begin{equation}
  \label{disq}
  S_{n} \; = \; \sum_{m=0}^{n} \, (-1)^{m} \, m! \, z^{-m}
\end{equation}
of the factorially divergent asymptotic series \eqref{ExpInt_AsySer} for
the exponential integral $E_{1} (z)$. We present the summation results
obtained for $z=3$ and $z=1/2$, respectively, by our new algorithm
\eqref{seps}, Wynn's epsilon algorithm \eqref{eps_al}, Brezinski's theta
algorithm \eqref{ThetaAl}, and the iteration \eqref{ThetIt} of
$\vartheta_{2}^{(n)}$.

\begin{table}
\caption{\label{divsq1}Summation of the asymptotic series ${}_{2} F_0 (1, 1; -1/z)=ze^zE_1(z)$ for $z=3$}
\begin{tabular}{|c|c|c|c|c|c|}\hline
$n$ & $\bigl\vert S_n-S \bigr\vert$ & $\bigl\vert
\tilde{T}_{2\Ent{(n-1)/2}}^{(n-1-2\Ent{(n-1)/2})}-S \bigr\vert$ & $\bigl\vert
\varepsilon_{2\Ent{n/2}}^{(n-2\Ent{n/2})}-S \bigr\vert$ & $\bigl\vert
\vartheta_{2\Ent{n/3}}^{(n-3\Ent{n/3})}-S \bigr\vert$ & $\bigl\vert
\mathcal{J}_{\Ent{n/3}}^{(n-3\Ent{n/3})}-S \bigr\vert$ \\
& \cref{disq}, z=3 & \cref{seps} & \cref{eps_al} &
\cref{ThetaAl} & \cref{ThetIt} \\ \hline
14 & 1.504$\times10^{4}$ & 2.587$\times10^{-8}$   & 9.326$\times10^{-7}$  &  5.357$\times10^{-11}$ & 4.735$\times10^{-11}$ \\
15 & 7.609$\times10^{4}$ & 3.287$\times10^{-9}$   & 7.006$\times10^{-7}$  &  1.379$\times10^{-11}$ & 6.865$\times10^{-12}$ \\
16 & 4.100$\times10^{5}$ & 1.568$\times10^{-9}$   & 2.877$\times10^{-7}$  &  4.047$\times10^{-11}$ & 6.424$\times10^{-11}$ \\
17 & 2.344$\times10^{6}$ & 3.844$\times10^{-11}$  & 2.192$\times10^{-7}$  &  3.737$\times10^{-10}$ & 3.744$\times10^{-11}$ \\
18 & 1.418$\times10^{7}$ & 8.860$\times10^{-9}$   & 9.445$\times10^{-8}$  &  4.069$\times10^{-11}$ & 1.258$\times10^{-12}$ \\
19 & 9.048$\times10^{7}$ & 2.204$\times10^{-9}$   & 7.289$\times10^{-8}$  &  6.707$\times10^{-11}$ & 2.000$\times10^{-10}$ \\
20 & 6.073$\times10^{8}$ & 4.265$\times10^{-9}$   & 3.272$\times10^{-8}$  &  3.673$\times10^{-13}$ & 1.386$\times10^{-12}$ \\
21 & 4.277$\times10^{9}$ & 2.101$\times10^{-10}$  & 2.552$\times10^{-8}$  &  9.726$\times10^{-14}$ & 1.785$\times10^{-12}$ \\
\hline
\end{tabular}
\end{table}

\begin{table}
\caption{\label{divsq2}Summation of the asymptotic series ${}_{2} F_0 (1, 1; -1/z)=ze^zE_1(z)$ for $z=1/2$}
\begin{tabular}{|c|c|c|c|c|c|}\hline

$n$ & $\bigl\vert S_n-S \bigr\vert$ & $\bigl\vert
\tilde{T}_{2\Ent{(n-1)/2}}^{(n-1-2\Ent{(n-1)/2})}-S \bigr\vert$ & $\bigl\vert
\varepsilon_{2\Ent{n/2}}^{(n-2\Ent{n/2})}-S \bigr\vert$ & $\bigl\vert
\vartheta_{2\Ent{n/3}}^{(n-3\Ent{n/3})}-S \bigr\vert$ & $\bigl\vert
\mathcal{J}_{\Ent{n/3}}^{(n-3\Ent{n/3})}-S \bigr\vert$ \\
&\cref{disq}, z=1/2 & \cref{seps} & \cref{eps_al} &
\cref{ThetaAl} & \cref{ThetIt}\\ \hline
14 & 1.379$\times10^{15}$ & 3.012$\times10^{-4}$   & 2.065$\times10^{-3}$  & 3.963$\times10^{-5}$   & 1.862$\times10^{-5}$ \\
15 & 4.147$\times10^{16}$ & 4.489$\times10^{-5}$   & 3.346$\times10^{-3}$  & 4.578$\times10^{-7}$   & 1.052$\times10^{-6}$ \\
16 & 1.330$\times10^{18}$ & 1.498$\times10^{-5}$   & 1.257$\times10^{-3}$  & 1.314$\times10^{-7}$   & 4.844$\times10^{-7}$ \\
17 & 4.529$\times10^{19}$ & 1.364$\times10^{-6}$   & 1.987$\times10^{-3}$  & 4.963$\times10^{-7}$   & 2.710$\times10^{-7}$ \\
18 & 1.633$\times10^{21}$ & 7.763$\times10^{-6}$   & 7.876$\times10^{-4}$  & 4.044$\times10^{-7}$   & 2.207$\times10^{-6}$ \\
19 & 6.214$\times10^{22}$ & 6.418$\times10^{-6}$   & 1.218$\times10^{-3}$  & 1.425$\times10^{-7}$   & 4.204$\times10^{-7}$ \\
20 & 2.489$\times10^{24}$ & 2.846$\times10^{-6}$   & 5.052$\times10^{-4}$  & 8.537$\times10^{-8}$   & 2.591$\times10^{-7}$ \\
21 & 1.047$\times10^{26}$ & 3.780$\times10^{-6}$   & 7.662$\times10^{-4}$  & 3.127$\times10^{-8}$   & 3.538$\times10^{-7}$ \\
\hline
\end{tabular}
\end{table}

%
\typeout{==> Section: Conclusions and discussions}
\section{Conclusions and discussions}
\label{Sec:ConclDiscuss}
%

Starting from the known recursion relations of Wynn's epsilon and rho
algorithm or of Osada's generalized rho algorithm, we constructed a
new convergence acceleration algorithm \eqref{seps} which is -- as shown
in our numerical examples in \cref{Sec:NumEx} -- not only successful in
the case of linearly convergent sequences, but works also in the case of
many logarithmically convergent sequences. Our numerical results also
showed that our new transformations sums factorially divergent
series. Thus, our new transformation \eqref{seps} is clearly more
versatile than either epsilon or rho, having similar properties as
Brezinski's theta algorithm \eqref{ThetaAl} or the iteration
\eqref{ThetIt} of $\vartheta_{2}^{(n)}$.

For the derivation of our new algorithm \eqref{seps} in
\cref{Sec:DiffEqTransOrd}, we pursued a novel approach. We analyzed the
finite difference equations with respect to the transformation order,
which the known recurrence formulas of Wynn's epsilon and rho algorithms
or of Osada's generalized rho algorithm satisfy. This approach yielded a
generalized recursive scheme \eqref{Rec_Tkn_Fn0_Fn1} which contains the
known recursions of Wynn's epsilon and rho algorithms and Osada's rho
algorithm as special cases.

To complete the derivation of our new algorithm \eqref{seps}, we only had
to choose appropriate initial conditions. Parts of the initial conditions
were chosen as in Wynn's epsilon algorithm according to
\cref{EpsIniVa}. In addition to these epsilon-type initial conditions
\eqref{EpsIniVa}, we then needed only one further initial condition for
$T_{2}^{(n)}$, which we identified with Brezinski's theta transformation
$\vartheta_{2}^{(n)}$ according to \cref{IniConds}. Thus, our new
algorithm \eqref{seps} can be viewed to be essentially a blend of Wynn's
epsilon and Brezinski's theta algorithm.

Of course, it should be possible to construct -- starting from the
generalized recursive scheme \eqref{Rec_Tkn_Fn0_Fn1} -- alternative new
sequence transformations by modifying our choice for $T_{2}^{(n)}$. The
feasibility of this approach will be investigated elsewhere.

\typeout{==> Section: Acknowledgements}
\section*{Acknowledgements}
\label{Sec:Ack}

Y. He was supported by the National Natural Science Foundation of China
(Grants no. 11571358), the China Postdoctoral Science Foundation funded
project (Grants no.\ 2012M510186 and 2013T60761), and the Youth
Innovation Promotion Association CAS. X.K. Chang was supported by the
National Natural Science Foundation of China (Grants no. 11701550,
11731014), and the Youth Innovation Promotion Association CAS. X.B. Hu
was supported by the National Natural Science Foundation of China (Grants
no. 11331008, 11571358).  J.Q. Sun was supported by the National Natural
Science Foundation of China (Grants no. 11401546, 11871444), and the
Shandong Provincial Natural Science Foundation, China (Grant no.\
ZR2013AQ024).

\newpage

\appendix
\appendixpage
\addappheadtotoc

\begin{appendix}

\typeout{==> Appendix: Reduced Bessel functions}

\section{Reduced Bessel functions}
\label{App:RBF}

In this Appendix, we briefly discuss the most important properties of the
so-called reduced Bessel functions \citep*[Eqs.\ (3.1) and
(3.2)]{Steinborn/Filter/1975c} and their anisotropic generalizations, the
so-called $B$ functions \citep*[Eqs.\ (3.3) and
(3.4)]{Filter/Steinborn/1978b}. These functions have gained some
prominence in molecular electronic structure theory (see for example
\cite{Weniger/2009a} and references therein), but have been largely
ignored in mathematics. Reduced Bessel and $B$ functions had also been
the topic of Weniger's Diploma and PhD theses
\cite{Weniger/1977,Weniger/1982}. But in this article, we are
predominantly interested in reduced Bessel functions because of the
series expansion \eqref{1/z_RBF} of $1/z$ in terms of reduced Bessel
functions, which is well suited to test the ability of a sequence
transformation to accelerate the logarithmic convergence of model
sequences of the type of \cref{ModSeqAlpha} with a non-integral decay
parameter.

Based on previous work by \citet[Eq.\ (55) on p.\ 15]{Shavitt/1963}, the
reduced Bessel function with in general complex order $\nu$ and complex
argument $z$ was introduced by \citet*[Eqs.\ (3.1) and
(3.2)]{Steinborn/Filter/1975c}:
\begin{equation}
  \label{Def:RBF}
  \widehat{k}_{\nu} (z) \; = \;
   (2/\pi)^{1/2} \, z^{\nu} \, K_{\nu} (z) \, .
\end{equation}
Here, $K_{\nu} (z)$ is a modified Bessel functions of the second kind
\citep*[Eqs.\ (10.27.4) and (10.27.5)]{Olver/Lozier/Boisvert/Clark/2010}.

Many properties of the reduced Bessel functions follow directly from the
corresponding properties of $K_{\nu} (z)$. One example is the recursion
\begin{equation}
  \label{Rec_RBF}
  \widehat{k}_{\nu+1} (z) \; = \; 2 \nu \, \widehat{k}_{\nu} (z) \, + \,
  z^{2} \, \widehat{k}_{\nu-1} (z) \, ,
\end{equation}
which is stable upwards and which follows directly from the three-term
recursion $K_{\nu+1} (z) = (2\nu/z) K_{\nu} (z) + K_{\nu-1} (z)$
\citep*[Eq.\ (10.29.1)]{Olver/Lozier/Boisvert/Clark/2010}. Another
example is the symmetry relationship
$\widehat{k}_{-\nu} (z) = z^{-2 \nu} \widehat{k}_{\nu} (z)$, which
follows from the symmetry relationship $K_{-\nu} (z) = K_{\nu} (z)$
\citep*[Eq.\ (10.27.3)]{Olver/Lozier/Boisvert/Clark/2010}.

If the order $\nu$ of $\widehat{k}_{\nu} (z)$ is half-integral,
$\nu = n + 1/2$ with $n \in \mathbb{N}_0$, a reduced Bessel function can
be expressed as an exponential multiplied by a terminating confluent
hypergeometric series ${}_1 F_1$ (see for example \cite[Eq.\
(3.7)]{Weniger/Steinborn/1983b}):
\begin{equation}
  \label{RBF_HalfInt}
  \widehat{k}_{n+1/2} (z) \; = \; 2^n \, (1/2)_n \,
   \mathrm{e}^{-z} \, {}_1 F_1 (-n; -2n; 2z) \, ,
    \qquad n \in \mathbb{N}_{0} \, .
\end{equation}
The functions $\widehat{k}_{n+1/2} (z)$ can be computed conveniently and
reliably with the help of the recursion (\ref{Rec_RBF}) in combination
with the starting values $\widehat{k}_{-1/2} (z) = \mathrm{e}^{-z}/z$ and
$\widehat{k}_{1/2} (z) = \mathrm{e}^{-z}$.

The polynomial part of $\widehat{k}_{n+1/2} (z)$ had been considered
independently in the mathematical literature. There, the notation
\begin{equation}
  \label{Def_BesPol}
  \theta_n (z) \; = \; \mathrm{e}^z \, \widehat{k}_{n+1/2} (z) \; = \;
  2^n \, (1/2)_n \, \mathrm{e}^{-z} \, {}_1 F_1 (-n; -2n; 2z)
\end{equation}
is used \cite[Eq.\ (1) on p.\ 34]{Grosswald/1978}. Together with some
other, closely related polynomials, the $\theta_n (z)$ are called
\emph{Bessel polynomials} \cite{Grosswald/1978}. According to
\citet[Section 14]{Grosswald/1978}, they have been applied in such
diverse and seemingly unrelated fields like number theory, statistics,
and the analysis of complex electrical networks.

Bessel polynomials occur also in the theory of Pad\'{e}
approximants. \citet[p.\ 82]{Pade/1892} had shown in his thesis that the
Pad\'{e} approximant $[n/m]_{\exp} (z)$ to the exponential function
$\exp (z)$ can be expressed in closed form, and
\citeauthor{Baker/Graves-Morris/1996} showed that Pad\'{e}'s expression
corresponds in modern mathematical notation to the following ratio of two
terminating confluent hypergeometric series ${}_1 F_{1}$ \citet[Eq.\
(2.12)]{Baker/Graves-Morris/1996}:
\begin{equation}
  [n/m]_{\exp} (z) \; = \; \frac
   {{}_1 F_{1} (- n; - n - m; z)}
   {{}_1 F_{1} (- m; - n - m; - z)} \, ,
   \qquad m, n \in \mathbb{N}_{0} \, .
\end{equation}
Thus, the \emph{diagonal} Pad\'{e} approximants $[n/n]_{\exp} (z)$ are
ratios of two Bessel polynomials \cite[Eq.\ (14.3-15)]{Weniger/1989}:
\begin{equation}
  [n/n]_{\exp} (z) \; = \; \frac{\theta_n (z/2)}{\theta_n (-z/2)} \, ,
   \qquad n \in \mathbb{N}_0 \, .
\end{equation}

The known monotonicity properties of the modified Bessel function
$K_{\nu} (z)$ imply that the reduced Bessel functions
${\widehat k}_{\nu} (z)$ are for $\nu > 0$ and $z \ge 0$ positive and
bounded by their values at the origin \cite[Eq.\
(3.1)]{Weniger/Steinborn/1983b}. In the case of reduced Bessel functions
with half-integral orders, this yields the following bound:
\begin{equation}
  \label{Boundedness_RBF}
  0 \, < \, \widehat{k}_{n+1/2} (z) \, \le \,
  \widehat{k}_{n+1/2} (0) \; = \; 2^n \, (1/2)_n \, ,
  \qquad 0 \le z < \infty \, . \quad n \in \mathbb{N}_{0} \, .
\end{equation}
In Grosswald's book \citep{Grosswald/1978}, it is shown that for fixed
and finite argument $z$ the Bessel polynomials $\theta_n (z)$ satisfy the
leading order asymptotic estimate \citep[p.\ 125]{Grosswald/1978}
\begin{equation}
  \label{Asy_BesPolTheta}
  \theta_n (z) \; \sim \; \frac {(2 n)!} {2^n n!} \, \mathrm{e}^{z}
   \; = \; 2^n \, (1/2)_n \, \mathrm{e}^{z} \, , \qquad n \to \infty \, .
\end{equation}
Combination of \cref{Def_BesPol,Asy_BesPolTheta} shows that the dominant
term of the Poincar\'{e}-type asymptotic expansion of
$\widehat{k}_{n+1/2} (z)$ with fixed and finite argument $z$ corresponds
to its value at the origin \citep[Eq.\ (3.9)]{Weniger/Steinborn/1983b}:
\begin{align}
  \label{AsyInd_RBF}
  \widehat{k}_{n+1/2} (z) & \; = \; \widehat{k}_{n+1/2} (0) \,
   \Bigl[1 + \mathrm{O} \bigl(n^{- 1} \bigr) \Bigr]
    \; = \; 2^n \, (1/2)_n \,
     \Bigl[1 + \mathrm{O} \bigl(n^{- 1} \bigr) \Bigr] \, ,
      \qquad n \to \infty \, .
\end{align}

For several functions, finite or infinite expansions in terms of reduced
Bessel functions are known. As a challenging problem for sequence
transformations, we propose to use the series expansion \citep[Eq.\
(6.5)]{Filter/Steinborn/1978a}
\begin{equation}
\label{1/z_RBF}
 \frac{1}{z} \; = \; \sum_{m=0}^{\infty} \,
  \widehat{k}_{m-1/2} (z) / [2^m \, m!] \, ,  \qquad z > 0 \, .
\end{equation}
\Cref{AsyInd_RBF} implies that the terms
$\widehat{k}_{m-1/2} (z) / [2^m m!]$ of this series decay like
$m^{- 3/2}$ as $m \to \infty$ \cite[p.\
3709]{Grotendorst/Weniger/Steinborn/1986}, or that the remainders
$\sum_{m=n+1}^{\infty} \, \widehat{k}_{m-1/2} (z) / [2^m m!]$ decay like
$n^{- 1/2}$ as $n \to \infty$. Thus, the series \eqref{1/z_RBF} converges
as slowly as the Dirichlet series
$\zeta (s) = \sum_{m=0}^{\infty} (m+1)^{s}$ with $s=1/2$, which is
notorious for extremely poor convergence. The slow convergence of the
infinite series \eqref{1/z_RBF} was demonstrated in \cite[Table
I]{Grotendorst/Weniger/Steinborn/1986}. After $10^{6}$ terms, only 3
decimal digits were correct, which is in agreement with the truncation
error estimate given above.

%
\typeout{==> Appendix: Levin's transformation}
\section{Levin's transformation}
\label{App:LevTr}
%

%
\typeout{==> Subsection: The general Levin transformation}
\subsection{The general Levin transformation}
\label{Sub:GenLevTr}
%

A sequence transformations of the type of Wynn's epsilon algorithm
(\ref{eps_al}) only requires the input of a finite sub-string of a
sequence $\{ S_{n} \}_{n=0}^{\infty}$. No additional information is needed for the
computation of approximations to the (generalized) limit $S$ of an input
sequence $\{ S_{n} \}_{n=0}^{\infty}$. Obviously, this is a highly advantageous
feature. However, in more fortunate cases, additional information on the
index dependence of the truncation errors $R_{n} = S_{n} - S$ is
available. For example, the truncation errors of Stieltjes series are
bounded in magnitude by the first terms neglected in the partial sums
(see for example \citep[Theorem 13-2]{Weniger/1989}), and they also have
the same sign pattern as the first terms neglected. The utilization of
such a \emph{structural} information in a transformation process should
enhance its efficiency. Unfortunately, there is no obvious way of
incorporating such an information into Wynn's recursive epsilon algorithm
(\ref{eps_al}) or into other sequence transformations with similar
features.

In \citeyear{Levin/1973}, \citet{Levin/1973} introduced a sequence
transformation which overcame these limitations and which now bears his
name. It uses as input data not only the elements of the sequence
$\{ S_{n} \}_{n=0}^{\infty}$, which is to be transformed, but also the
elements of another sequence $\{ \omega_{n} \}_{n=0}^{\infty}$ of
explicit \emph{estimates} of the remainders $R_{n} = S_{n} - S$. These
remainders estimates, which must be explicitly known, make it possible to
incorporate additional information into the transformation process and
are thus ultimately responsible for the remarkable versatility of Levin's
transformation.

For Levin's transformation, we use the notation introduced in \citep[Eqs.\
(7.1-6) and (7.1-7)]{Weniger/1989}:
\begin{align}
  \label{LevTrDiffOpRep}
  \mathcal{L}_{k}^{(n)} (\beta, S_{n}, \omega_{n}) & \; = \;
   \frac {\Delta^{k} [(\beta+n)^{k-1} S_{n} / \omega_{n}]}
    {\Delta^{k} [(\beta+n)^{k-1} / \omega_{n}]}
  \\
  \label{GenLevTr}
  & \; = \; \frac {\displaystyle \sum_{j=0}^{k} \, (-1)^{j} \,
   {\binom{k}{j}} \, \frac {(\beta+n+j)^{k-1}} {(\beta+n+k)^{k-1}} \,
    \frac {S_{n+j}} {\omega_{n+j}}} {\displaystyle \sum_{j=0}^{k} \,
     (-1)^{j} \, {\binom{k}{j}} \, \frac {(\beta+n+j)^{k-1}}
      {(\beta+n+k)^{k-1}} \, \frac {1} {\omega_{n+j}} } \, ,
       \qquad k, n \in \mathbb{N}_{0} \, , \quad \beta > 0 \, .
\end{align}
Here, $\beta > 0$ is a shift parameter, $\{ S_{n} \}_{n=0}^{\infty}$ is the input
sequence, and $\{ \omega_{n} \}_{n=0}^{\infty}$ is a sequence of remainder estimates.

In \cref{LevTrDiffOpRep,GenLevTr} it is essential that the input sequence
$\{ S_{n} \}$ starts with the sequence element $S_{0}$. The reason is
that both the definitions according to \cref{LevTrDiffOpRep,GenLevTr} as
well as the recurrence formulas for its numerator and denominator sums
(see for example\citep[Eq.\ (3.11)]{Weniger/2004}) depend explicitly on
$n$ as well as on $\beta$.

Levin's transformation, which is also discussed in the NIST Handbook
\citep*[\S 3.9(v) Levin's and Weniger's
Transformations]{Olver/Lozier/Boisvert/Clark/2010}, is generally
considered to be both a very powerful and a very versatile sequence
transformation (see for example \citep*{Borghi/Weniger/2015,%
  Brezinski/RedivoZaglia/1991a,Sidi/2003,Smith/Ford/1979,%
  Smith/Ford/1982,Weniger/1989,Weniger/2004} and references therein). The
undeniable success of Levin's transformation inspired others to construct
alternative sequence transformations that also use explicit remainder
estimates (see for example
\citep*{Cizek/Zamastil/Skala/2003,Homeier/2000a,Weniger/1989,%
  Weniger/1992,Weniger/2004} and references therein).

We still have to discuss the choice of the so far unspecified remainder
estimates $\{ \omega_{n} \}_{n=0}^{\infty}$. A principal approach would be to look for
remainder estimates that reproduce the leading order asymptotics of the
actual remainders \citep[Eq.\ (7.3-1)]{Weniger/1989}:
\begin{equation}
  \label{AsyCondRemEst}
  R_{n} \; = \; S_{n} \, - \, S \; = \; \omega_{n} \,
   \left[ C + \mathrm{O} \bigl( 1/n \bigr) \right] \, ,
    \qquad C \ne 0 \, , \quad n \to \infty \, .
\end{equation}
This is a completely valid approach, but there is the undeniable problem
that for every input sequence $\{ S_{n} \}_{n=0}^{\infty}$ the leading order asymptotics
of the corresponding remainders $R_{n} = S_{n} - S$ as $n \to \infty$ has
to be determined. Unfortunately, such an asymptotic analysis may well
lead to some difficult technical problem and be a non-trivial research
problem in its own right.

In practice, it is much more convenient to use simple explicit remainder
estimates introduced by \citet{Levin/1973} and \citet*{Smith/Ford/1979},
respectively, which are known to work well even in the case of purely
numerical input data. In this article, we only consider Levin's remainder
estimates $\omega_{n} = (\beta + n) \Delta S_{n-1}$, and
$\omega_{n} = - [\Delta S_{n-1}][\Delta S_{n}]/[\Delta^{2} S_{n-1}]$,
which lead to Levin's $u$ and $v$ variants \citep[Eqs.\ (58) and
(67)]{Levin/1973}. These estimates lead to transformations that are known
to be powerful accelerators for both linear and logarithmic convergence
(compare also the asymptotic estimates in \citep[Section
IV]{Weniger/2004}). Therefore, Levin's $u$ and $v$ transformations are
principally suited to represent the remaining initial condition
$T_{2}^{(n)}$ in \cref{Choice_F0F1}.

%
\typeout{==> Subsection: Levin's $u$ transformation}
\subsection{Levin's \protect{$u$} transformation}
\label{Sub:LevinU}
%

Levin's $u$ transformation, which first appeared already in
\citeyear{Bickley/Miller/1936} in the work of
\citet{Bickley/Miller/1936}, is characterized by the remainder estimate
$\omega_{n} = (\beta + n) \Delta S_{n-1}$ \citep[Eq.\ (58)]{Levin/1973}.
Inserting this into the finite difference representation
\eqref{LevTrDiffOpRep} yields:
\begin{equation}
  \label{uLevTrDiffOpRep}
  u_{k}^{(n)} (\beta, S_{n}) \; = \;
   \frac {\Delta^{k} [(\beta+n)^{k-2} S_{n} / \Delta S_{n-1}]}
    {\Delta^{k} [(\beta+n)^{k-2} / \Delta S_{n-1}]} \, .
\end{equation}
For $k=2$, the right-hand side of this expression depends on $n$ only
\emph{implicitly} via the input data $\{ S_{n} \}$ and also does not
depend on the scaling parameter $\beta$. Thus, we obtain:
\begin{equation}
  \label{uLevTr_2}
  u_{2}^{(n)} (\beta, S_{n}) \; = \;
   \frac {\Delta^{2} \bigl[ S_{n} / \Delta S_{n-1} \bigr]}
    {\Delta^{2} \bigl[ 1 / \Delta S_{n-1} \bigr]} \, .
\end{equation}
If we now use $S_{n} = S_{n-1} + \Delta S_{n-1}$, we obtain
$S_{n} / \Delta S_{n-1} = S_{n-1} / \Delta S_{n-1} + 1$, which implies
$\Delta [S_{n}/\Delta S_{n-1}] = \Delta [S_{n-1}/\Delta S_{n-1}]$. Thus,
we obtain the alternative expression,
\begin{equation}
  \label{uLevTr_3}
  u_{2}^{(n)} (\beta, S_{n}) \; = \;
   \frac {\Delta^{2} [ S_{n-1} / \Delta S_{n-1} ]}
    {\Delta^{2} [ 1 / \Delta S_{n-1} ]} \, , 
\end{equation}
which can be reformulated as follows:
\begin{equation}
  \label{uLevTr_4}
  u_{2}^{(n)} (\beta, S_{n}) \; = \; \frac
   {S_{n} \bigl[ \Delta S_{n+1} \bigr] \bigl[ \Delta^{2}  S_{n-1} \bigr]
    -  S_{n+1} \bigl[ \Delta S_{n-1} \bigr]
     \bigl[ \Delta^{2} S_{n} \bigr]}
   {\bigl[ \Delta S_{n+1} \bigr] \bigl[ \Delta^{2} S_{n-1} \bigr] -
    \bigl[ \Delta S_{n-1} \bigr] \bigl[ \Delta^{2} S_{n} \bigr]} \, .
\end{equation}
In order to enhance numerical stability, it is recommendable to rewrite
this expression as follows:
\begin{equation}
  \label{uLevTr_5}
  u_{2}^{(n)} (\beta, S_{n}) \; = \; S_{n} \, - \, \frac
   {\bigl[ \Delta S_{n-1} \bigr] \bigl[ \Delta S_{n} \bigr]
    \bigl[ \Delta^{2} S_{n} \bigr]}
   {\bigl[ \Delta S_{n+1} \bigr] \bigl[ \Delta^{2} S_{n-1} \bigr] -
    \bigl[ \Delta S_{n-1} \bigr] \bigl[ \Delta^{2} S_{n} \bigr]} \, .
\end{equation}

%
\typeout{==> Subsection: Levin's $v$ transformation}
\subsection{Levin's \protect{$v$} transformation}
\label{Sub:LevinV}
%

Levin's $v$ transformation is characterized by the remainder estimate
\citep[Eq.\ (67)]{Levin/1973} which is inspired by Aitken's $\Delta^{2}$
formula \eqref{AitFor_1}:
\begin{equation}
  \label{vRemEst}
  \omega_{n} \; = \; \frac{[\Delta S_{n-1}] [\Delta S_{n}]}
   {\Delta S_{n-1} - \Delta S_{n}} \; = \;
    - \frac{[\Delta S_{n-1}] [\Delta S_{n}]}{\Delta^{2} S_{n-1}}
    \; = \; \frac{[\Delta S_{n-1}] [\Delta S_{n}]}{\Delta^{2} S_{n-1}} \, ,
\end{equation}
In \cref{vRemEst}, we made use of the fact that Levin's transformation
$\mathcal{L}_{k}^{(n)} (\beta, S_{n}, \omega_{n})$ is a homogeneous
function of degree zero of the $k+1$ remainder estimates $\omega_{n}$,
$\omega_{n+1}$, \dots. Thus, we can use any of the expressions in
\cref{vRemEst} without affecting the value of the transformation.

Inserting the remainder estimate \eqref{vRemEst} into the finite
difference representation \eqref{LevTrDiffOpRep} yields:
\begin{equation}
  \label{vLevTrDiffOpRep}
  v_{k}^{(n)} (\beta, S_{n}) \; = \; \frac
   {\displaystyle \Delta^{k} \frac{\displaystyle
    (\beta+n)^{k-1} S_{n} [\Delta^{2} S_{n-1}]}{\displaystyle [\Delta
     S_{n-1}] [\Delta S_{n}]}}
   {\displaystyle \Delta^{k} \frac{\displaystyle
    (\beta+n)^{k-1} [\Delta^{2} S_{n-1}]}{\displaystyle [\Delta
     S_{n-1}] [\Delta S_{n}]}} \, .
\end{equation}
For $k=1$, the right-hand side of this expression depends on $n$ only
\emph{implicitly} via the input data $\{ S_{n} \}$ and also does not
depend on the scaling parameter $\beta$. We obtain:
\begin{equation}
  \label{vLevTr_1}
  v_{1}^{(n)} (\beta, S_{n}) \; = \; \frac
   {\displaystyle  S_{n+1} [\Delta S_{n-1}] [\Delta^{2} S_{n}] - S_{n}
    [\Delta S_{n+1}] [\Delta^{2} S_{n-1}]}
    {\displaystyle [\Delta S_{n-1}] [\Delta^{2} S_{n}] -
    [\Delta S_{n+1}] [\Delta^{2} S_{n-1}]} \, .
\end{equation}
Comparison of \cref{uLevTr_4,vLevTr_1} yields
$v_{1}^{(n)} (\beta, S_{n}) = u_{2}^{(n)} (\beta, S_{n})$ for arbitrary
$\beta \in \mathbb{R}$ and $n \in \mathbb{N}_{0}$.

Again, it is recommendable to rewrite \cref{vLevTr_1} as follows:
\begin{equation}
  \label{vLevTr_2}
  v_{1}^{(n)} (\beta, S_{n}) \; = \; S_{n} + \frac
   {\displaystyle [\Delta S_{n}] [\Delta S_{n-1}] [\Delta^{2} S_{n}]}
    {\displaystyle [\Delta S_{n-1}] [\Delta^{2} S_{n}] -
     [\Delta S_{n+1}] [\Delta^{2} S_{n-1}]} \, .
\end{equation}

\end{appendix}

\addcontentsline{toc}{section}{Bibliography}

\end{document}